\documentclass[11pt,notcite,notref]{article}

\usepackage[T1]{fontenc}
\usepackage[utf8]{inputenc}
\usepackage{amsmath,cases}
\usepackage{amsfonts, amssymb}
\usepackage{graphicx}%
\usepackage{placeins}
\usepackage{indentfirst}
\usepackage{caption}
\usepackage{subcaption}
\usepackage{bm}
\usepackage{xcolor}
\usepackage{tikz}
\usetikzlibrary{decorations.pathreplacing,decorations.markings}
\usepackage{comment}
\usepackage{cancel}
\usepackage{empheq}
\allowdisplaybreaks[1]
\usepackage[a4paper, total={6in, 8in}]{geometry}

\numberwithin{equation}{section}

\providecommand{\U}[1]{\protect\rule{.1in}{.1in}}
\newtheorem{theorem}{Theorem}

\newtheorem{remark}[theorem]{Remark}

\renewcommand{\div}{\mathrm{div}_{\mathbf{x}}}
\newcommand{\vc}[1]{\mathbf{#1}}
\newcommand{\vm}[1]{\mathbb{#1}}

\newcommand{\ignore}[1]{}

\usepackage{authblk}
\AtBeginDocument{}

\date{}
\title{Revisiting Shikhmurzaev's approach to the contact line problem}
\author[1]{Amrita Ghosh \thanks{ghosh@iam.uni-bonn.de}}
\affil[1]{Institute for Applied Mathematics, University of Bonn \authorcr Endenicher Allee 60, 53115 Bonn, Germany}
\author[1]{Barbara Niethammer \thanks{niethammer@iam.uni-bonn.de}}
\author[1]{Juan J.L. Vel\'azquez \thanks{velazquez@iam.uni-bonn.de}}

\begin{document}
	\maketitle
	
	\begin{abstract}
		In this paper, we revisit a model for the contact line problem which has been proposed 
		by Shikhmurzaev (1993). 
		In the first part, in addition to rederiving the model, we study in detail the assumptions required to obtain the isothermal limit of the model. We also derive in this paper several lubrication approximation models, based on Shikhmurzaev's approach. The first two lubrication models describe thin film flow of incompressible fluids on solid substrates, 
		based on different orders of magnitude of the slip length parameter. The third lubrication model describes a meniscus formation where a wedge-shaped solid immerses in a thin film of fluid.
		
	\end{abstract}
	
	\section{Introduction}

	In this paper, we revisit Shikhmurzaev's theory for fluid interfaces \cite{shikh93}. In particular, we discuss in detail the conditions required to be able to assume that the temperature of the system is constant. We derive also some thin film approximations of the model.
	
To explain the background briefly, a 'contact line' is formed when three distinct phases of matter meet. If the phases are all fluid (liquid or gas) phases, that is, two immiscible fluids are in contact with another liquid or a gas, this contact line is freely deformable in space, while it is bound to move on a given surface, if one of the phases is a solid. We are interested in this latter case. The process of displacement of one fluid by another immiscible fluid from a solid surface in the framework of classical fluid mechanics is known as the 'moving contact line' phenomenon. This contact line motion appears in many natural and technological processes, for example, coating of solids by liquid films, photographic films or polymer processing.
	
	The mathematical modelling of the dynamic moving contact line encounters a fundamental difficulty, compared to the equilibrium state. 
	The typical no-slip boundary condition for velocities of both the fluids on the solid surface 
	gives rise to a non-integrable singularity of the shear stress, that is, to an infinite force exerted by the fluid on the solid (cf. \cite{HS}, \cite{DD}) which is considered not to be physical. The pressure of the fluid at the contact line becomes singular as well.
	This problem is known in the literature as the "moving contact line paradox". 
	
	Being important from both theoretical and practical point of view, this problem has been investigated from different perspectives. In order to have a well-defined problem in the presence of a contact line, one must address two fundamental questions: (1) how to determine the (dynamic) contact angle i.e. the angle between the free surface and the solid boundary as the contact line moves, and (2) how to remove the divergence of the dissipation of the energy formula. In determining the contact angle, it is important to take into account that the contact angle is not merely a material property, but it depends also on the velocity of the contact line (as suggested by experiments).
	The survey of Dussan \cite{dussan} discusses different approaches by physicists and chemists in determining the dynamics of a contact line. It is broadly accepted that the dynamic and equilibrium contact angles $\theta_d$ and $\theta_e$ are related via
	\begin{equation}
		\label{dynamic_CA}
		\vc{v}^c = F(\cos (\theta_d) - \cos(\theta_e)), 
	\end{equation}
	where $\vc{v}^c$ is the contact line velocity (along the solid), and $F$ is some increasing function such that $F(0)=0$.
	But there have been different suggestions for the precise form of $F$ in equation (\ref{dynamic_CA}). Blake-Haynes \cite{BH} derived (\ref{dynamic_CA}) by using thermodynamic and molecular kinetics while Cox \cite{cox'86} used matched asymptotic analysis and hydrodynamic arguments. Also Ren-E \cite{RenE11} derived (\ref{dynamic_CA}) from thermodynamic principles applied to constitutive equations.
	On the other hand,
	to remove the shear-stress singularity near the contact line, a widely used approach is to consider the Navier boundary condition which says that the slip of the fluid along the solid is proportional to the tangential stress acting from the fluid on the fluid-solid interface, i.e.
	\begin{equation*}
		\vc{v}^L \cdot \vc{n}=\vc{v}^S \cdot \vc{n}, \quad \tau \cdot \mathbb{S} \cdot \vc{n}= \beta \left( \vc{v}^L - \vc{v}^S\right) \cdot \tau,
	\end{equation*}
	with some parameter $\beta >0$ which is known as the \textit{friction parameter}, or the \textit{slip coefficient}. Variants of the above slip condition concerning different assumptions on $\beta$ have been considered in a number of works, for example in \cite{hocking77, huh_mason77, greenspan78, baiocchi90}. 
	The above references are by no means exhaustive and we refer to \cite[Section 3.4]{shikh97} for further reading.  
	
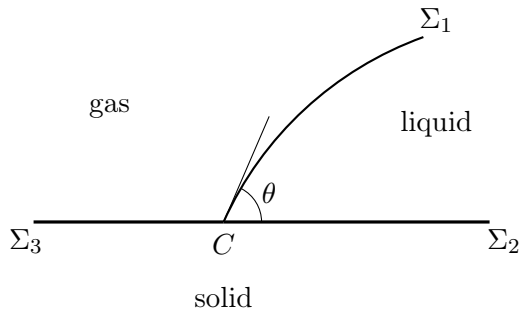
\begin{figure}[h]
		\centering
		\begin{tikzpicture}
			\draw[very thick] (-2.5,0) -- (0,0);
			\draw (0,0) -- (0.6,1.4);
			\draw[very thick,postaction={decorate}] (0,0) -- (3.5,0);
			\node at (0,-0.29) {$C$};
			\node at (2.8, 1.3) {liquid};
			\node at (-1.5, 1.5) {gas};
			\node at (0, -1) {solid};
			\node at (2.8, 2.68) {$\Sigma_1$};
			\node at (-2.6, -0.25) {$\Sigma_3$};
			\node at (3.69, -0.25) {$\Sigma_2$};
			\draw (0.5,0) arc (0:62:0.5);
			\node at (0.6, 0.39) {$\theta$};
			\draw[thick] (0,0) arc (156:110:4.6);
		\end{tikzpicture}
		\caption{General three-phase model with contact line $C$}
		\label{fig1}
	\end{figure}

	Another model has been introduced by Y. Shikhmurzaev in \cite{shikh93}, based on the approach proposed by Bedeaux et al \cite{BAM}, to generalize the methods of irreversible thermodynamics
	to systems containing interfaces.
	The basic idea of Shikhmurzaev's theory is to consider a three phase model (cf. Figure \ref{fig1}) and describe the interfaces as surfaces of zero thickness, separating the bulk phases. Although the thickness of these interfacial layers is negligible compared to the characteristic length scale of the bulk flow, they can possess some intrinsic surface properties (e.g. surface tension) which play an important role in the overall dynamics via capillary effects. This approach formulates in a systematic way the boundary conditions at the interfaces for the equations describing the bulk phases. Similarly, the contact line
	can be defined as boundaries of the interfaces.
	Incorporating the fundamental laws of fluid mechanics and the theory of irreversible thermodynamics, in addition to the classical conservation equations of mass, momentum and energy, different constitutive equations and transport coefficients are obtained, giving the required boundary conditions. In order to model the thermodynamics of the interfaces, the idea of Onsager's principle has been used, as for the bulk phases. 
	As a consequence, the (Navier) slip boundary condition is deduced at the fluid-solid interface which makes the force on the solid surface bounded, solving the paradoxical phenomena described in the beginning. Note that this derivation also makes clear that the slip boundary condition is not just an empirical assumption as often thought of, but rather can be formulated rigorously. Secondly, the dynamic contact angle is determined as part of the solution, given by a generalization of Young's law.
	
	The main novelty of this model is that it incorporates the aspect of interface formation and disappearance at the contact line, apart from addressing the two fundamental issues of contact line problem. This model is very appealing and has a great potential to explain several physical phenomena such as wetting, coalescence, cusp formation etc (cf. the introduction of \cite{bil06}).
	
	 In Section \ref{S2}, we rederive the model in situations which might contain contact lines at the boundary of the solid, liquid and gas phases. We made the derivation as mathematically rigorous as possible, using the physical ideas of Bedeaux et al \cite{BAM}, as well as Shikhmurzaev's. One of the issues that we try to clarify is the boundary conditions at the contact lines. It was indicated in \cite{bil06} that there is a missing boundary condition which is not included in the original derivation of Shikhmurzaev (in \cite{shikh93}, \cite{shikh97}). In \cite{bothe}, the authors also indicated the requirement of an additional boundary condition, although the one derived in \cite{bothe} is slightly different from the one in \cite{bil06}. In our derivation, we obtain the need for the additional boundary condition and we obtain the same one suggested in \cite{bothe}. This additional condition is required in order to obtain a well-posed mathematical problem.
	

	In Section \ref{S3}, we derive different lubrication models using Shikhmurzaev's approach to the contact line problem.
	Here the motivation is to obtain in the particular setting of thin films, a simpler reduced model than the original one, in which the differences between this approach and the other solution models of the contact line problem (for example, Navier boundary condition, disjoining pressure, precursor film etc) can be compared.
	Thin film equations have been discussed in many works, both experimentally and theoretically, taking into account different theories and approaches, and to answer different questions.
	Derivation of the classical thin film equation in case of a droplet can be found in, for example \cite{greenspan78, HS}, while in \cite{GGO}, the authors have considered thin film equations introducing the slip parameter in the model. On the other hand, there are other lubrication models, such as of a thin fluid in the Hele-Shaw cell \cite{CPG}, of Darcy flow \cite{GO}, rupture of thin film subject to van der Waals forces \cite{ED}, considering different orders of magnitude of the slip-length parameter \cite{MWW}, for an active gel \cite{KMW}, for surfactant driven viscous film \cite{GW}.
	In \cite{shikh2005}, the author has discussed a rupture phenomenon in the framework of the thin film approximation.
	Also an approximation of Shikhmurzaev's model for a gravity driven flow on an inclined plane has been considered in \cite{bil08}. In particular, the asymptotic behaviour of the solutions in a thin film setting has been computed.
We refer to the review paper \cite{W} for an overview.
	The classical work showing the existence of a weak solution corresponding to the zero contact angle goes back to \cite{BF90}. Also we mention the recent work \cite{GT} where the authors establish the existence of a global in time solution of the full contact line problem, provided the data is sufficiently close to the equilibrium value and the equilibrium contact angle may vary between $0$ and $\pi$.
	
	\section{Shikhmurzaev's model revisited}
	\label{S2}
	
	\subsection{Setting}
	In the following, we summarize the main ideas in the derivation of Shikhmurzaev's model. Specifically, we reformulate the model in the notations that we will use in rest of the paper and discuss in  detail some of the assumptions, in particular the isothermal limit, that has been used in most of Shikhmurzaev's work. Let us consider a simple model consisting of a liquid, gas and solid phase (more general situations with two liquids can be found in \cite{shikh97}). The model includes the evolution of interfaces separating two bulk phases.
	
	Let us consider a liquid with density $\varrho^L$.
	As we assume that our system can be described as a continuous medium, the properties of the fluid can be completely characterized by the velocity $\vc{v}^L$ and the temperature $T$ at each point $(x,t)$. Concerning the gas, it has density $\varrho^G$ and pressure $p^G$ at each point, while for the solid body, we assume that it is union of a finite number of rigid pieces and therefore we do not need to describe its deformations. Nevertheless, one could consider situations where there is a relative motion between the different parts that composes the solid boundary, for instance a container whose part of the boundary is a mobile piston. The rigid solid has a temperature $T$, as well as density $\varrho^S$. 
	
	The interaction between the different phases takes place through the liquid-gas, liquid-solid and gas-solid interfaces which are denoted by $\Sigma_1, \Sigma_2, \Sigma_3$ respectively. The three interfaces meet at a line $C$ which is termed as the contact line (cf. Figure \ref{fig1}). We assume in the following that all the interfaces $\Sigma_i$ are contained in smooth surfaces and that the contact line $C$ is a smooth curve. More precisely, we assume that the gas and the liquid fill a domain $\Omega$ with smooth boundary $\partial\Omega$
	. Then $\partial\Omega = \Sigma_2\cup \Sigma_3$. Note that, in general, the domain $\Omega$ changes in time. This change of $\Omega$ happens if the domain translates rigidly in space or if $\Omega$ changes its shape, or a combination of translation and deformation of $\Omega$. We can also have the situation where $\Omega$ does not change in time, but $\partial\Omega$ has a non-zero velocity. This happens for instance in the case in which the boundary of the container is a rotating cylinder. We characterize in the following the motion of the solid domain by a velocity function $\vc{v}^S$ defined on $\partial\Omega$.
	We will assume in the following, $\partial\Omega(t) = \{(x,t)\in\mathbb{R}^3\times \mathbb{R}, x=\Gamma(\sigma,t),\sigma\in\partial\Omega(0)\}$. The parameter $\sigma$ characterizes uniquely each material point of the boundary. Then the velocity of each boundary point is given by $\vc{v}^S(\sigma,t)= \frac{\partial\Gamma}{\partial t}(\sigma,t)$. Notice that the assumption that the boundary of the container consists of rigid parts imposes stringent constraints on $\Gamma(\sigma,t)$ and $\vc{v}^S(\sigma,t) $. We remark also that the domain $\Omega$ does not change in time if $\vc{v}^S\cdot \vc{n} =0$ where $\vc{n}$ is the normal vector to the boundary $\partial\Omega$. 
	
	These interfaces $\Sigma_i$ can be thought of as two-dimensional continuous media that can be described by a set of macroscopic variables. In our problem, the variables that we use to characterize these interfaces are the local surface densities $\varrho^s_i, i=1,2,3$, the velocities $\vc{v}^s_i, i=1,2,3$ and the temperature $T$. We emphasize that the surface densities $\varrho^s_i$ are not mass per unit volume, but mass per unit surface. Notice that we expect the mass for the interfaces to be small compared to the bulk, i.e. if $\ell$ denotes the thickness of the interface, then $\varrho^s_i\simeq \varrho^L \ell$. Moreover, denoting by $L$ the characteristic length scale of the full system, the mass in the liquid is given by $m^L = \varrho^L L^3$ while the mass of the interface is $m^s_i = \varrho^s_i L^2$. Therefore, $\frac{m^s_i}{m^L} = \frac{\varrho^s_i}{\varrho^L L} = \frac{\ell}{L} \ll 1$. However, the interfaces are able to produce relevant mechanical forces (for example, surface tension) that can be strongly dependent on the densities $\varrho^s_i$, which play an important role in the properties and dynamics of the interfaces. This is the reason that we keep $\varrho^s_i$ in the model. A relevant feature of Shikhmurzaev's approach is that $\varrho^s_i$ can take different values at different points of the interfaces.
	
	Concerning the contact line, the velocity is denoted by $\vc{v}^c$, defined as in (\ref{contact_line_velocity}).
	
	At this point, we assume that the system is in local equilibrium, i.e. all the thermodynamic relations (such as Gibbs relation) hold for every infinitesimal region of the space both in the bulk and at the interfaces.
	
	\subsubsection{Assumptions}
	\label{S.assump}
	
	The main assumptions that we will use in the derivation are the following.
	We restrict ourselves to the case of an incompressible liquid, i.e. the liquid density $\varrho^L$ is assumed to be constant (the compressible case could be considered similarly).
	
	Since the gas density is much lower than that of the liquid, we would like to consider $\varrho^G =0$.
	On the other hand, we want to have a gas with pressure comparable to that of the liquid. These gas properties have to be understood in a suitable asymptotic sense and the precise meaning and the consistency of these properties is explained in detail later in Section \ref{Conservation_Law}. We remark here that the assumption $\varrho^G \ll \varrho^L$ leads us to assume that both the internal and kinetic energy of the gas per unit of volume are negligible compared to the analogous quantities for the liquid. We will also see later 
	that the heat transport in the gas is much smaller than in the liquid for similar differences of temperature.
	
	The density of the solid body $\varrho^S$ is also assumed to be constant. As the solid is considered rigid, its velocity $\vc{v}^S$ can be characterized by the translational velocity $\vc{U}(t)$ and rotational velocity $\omega(t)$, thus, 
	\begin{equation}
		\label{solid_velocity}
		\vc{v}^S = \vc{U}(t) + \omega(t) \wedge \vc{x}.
	\end{equation}
	
	We also assume that the temperature across all the phases, the bulk, at the interface and at the contact line, is continuous, but not necessarily constant. That is, for $x_0\in \Sigma_i, i=1,2,3$,
	\begin{equation}
		\label{temperature_continuity}
		T(x_0) = \lim_{x\to x_0, x\in \Omega_L} T(x) = \lim_{x\to x_0, x\in \Omega_S} T(x),
	\end{equation}
	but $T(x_0)\ne T(x_1), x_0,x_1\in\Sigma_i$ in general.
	
	We further assume that the total mass of the contact line is negligible, and also that it contains neither momentum nor energy. However, suitable conditions concerning the conservation of mass, momentum and energy hold on it. This is not the only possibility; In the re-formulation of Shikhmurzaev's model in \cite{bothe}, the contact line density has been assumed to be non-zero.
	
	In the Appendix, we precise the notations used in this work and recall the necessary definitions and formulas from differential geometry. Throughout this article, we use the superscripts $L, S, G, s$ to denote the quantities in the liquid, solid, gas and at the interfaces respectively.
	
	\subsection{Conservation laws }
	\label{Conservation_Law}
	
	With the previous assumptions, the total mass density of the system becomes,
	\begin{equation}
		\label{mass}
		\varrho := \varrho^L \chi^L +\varrho^S \chi^S + \varrho^s_1 \delta_{\Sigma_1} +  \varrho^s_2 \delta_{\Sigma_2} +  \varrho^s_3 \delta_{\Sigma_3},
	\end{equation}
	where $\chi^L, \chi^S$ are the characteristic functions of the region filled by the liquid and the solid respectively and $\delta_{\Sigma_i}$ denotes Dirac measures, supported on the surfaces $\Sigma_i, i=1,2,3$. For later use, we also introduce the characteristic function $\chi^G$ of the region occupied by the gas. Notice that $\partial(\text{supp}\;\chi^L) = \Sigma_1\cup \Sigma_2$.
	
	Similarly, the total momentum density of the system is given by,
	\begin{equation}
		\label{momentum}
		\vc{M} := \varrho^L\vc{v}^L\chi^L + \varrho^S\vc{v}^S\chi^S + \varrho^s_1 \vc{v}^s_1 \delta_{\Sigma_1} + \varrho^s_2 \vc{v}^s_2 \delta_{\Sigma_2} + \varrho^s_3 \vc{v}^s_3 \delta_{\Sigma_3}.
	\end{equation}
	Note that in quantities (\ref{mass}) and (\ref{momentum}), the contribution of the gas to the mass and momentum density has been completely neglected, due to the above assumptions in Section \ref{S.assump}.
	
	Further, the total energy density is
	\begin{equation}
		\label{energy}
		E:= \varrho^L e^L \chi^L +\varrho^S e^S \chi^S 
		+ \varrho^s_1 e^s_1 \delta_{\Sigma_1}+ \varrho^s_2 e^s_2\delta_{\Sigma_2}+ \varrho^s_3 e^s_3\delta_{\Sigma_3}
	\end{equation}
	with
	$$
	e^\xi = \frac{1}{2}|\vc{v}^\xi|^2 + u^\xi, 
	\qquad e^s_i = \frac{1}{2}|\vc{v}^s_i|^2 + u^s_i, \quad \xi = L,S,G, \ i= 1,2,3.
	$$
	Here $e^\xi, e^s_i$ are the total energies per unit mass, $\frac{1}{2}|\vc{v}^\xi|^2, \frac{1}{2}|\vc{v}^s_i|^2$ are the kinetic energies per unit mass and $u^\xi, u^s_i$ are the internal energies per unit mass, for the bulk and the interfaces respectively.
	
	We neglect in (\ref{energy}) the internal energy of the gas compared with the one of the liquid. Indeed, for the changes of internal energy per unit mass, we have \begin{equation}
		\label{internal_energy_temperature}
		\delta u^G = c^G_m \;\delta T, \qquad \delta u^L = c^L_m \; \delta T
	\end{equation}
	where $c^G_m, c^L_m$ are the specific heat capacities of the gas and the liquid respectively. To obtain the above quantities for unit volume, we must multiply with the densities of the gas and the liquid respectively. Since the specific heat capacities of gas and liquid are of the same order of magnitude\footnote{ cf. https://en.wikipedia.org/wiki/Table\_of\_specific\_heat\_capacities}, it follows that $\varrho^G u^G \ll \varrho^L u^L$. Also, the kinetic energy per unit volume of the gas is much smaller compared to the liquid if $\vc{v}^G, \vc{v}^L$ are comparable, i.e. $\frac{1}{2}\varrho^G|\vc{v}^G|^2 \ll \frac{1}{2}\varrho^L|\vc{v}^L|^2 $. Therefore, we get that $\varrho^G e^G \ll \varrho^L e^L$ and thus, we ignore the contribution of the gas to the energy and hence to the entropy.
	
	We now derive the evolution equations for our system. These are obtained from the following basic conservation equations for mass, momentum and energy:
	\begin{align}
		\label{mass_conv}
		\partial_t \varrho + \div \;\vc{M} &=0,\\
		\label{momentum_conv}
		\partial_t \vc{\vc{M}} + \div \; \vc{J}_M& =0,\\
		\label{energy_conv}
		\partial_t E + \div \;\vc{J}_E&=0.
	\end{align}
	Here, the momentum flux is defined as,
	\begin{equation}
		\label{momentum_flux}
		\begin{aligned}
			\vc{J}_M := &\left( \varrho^L \vc{v}^L \otimes \vc{v}^L- \mathbb{S}^L\right) \chi^L +\left( \varrho^S \vc{v}^S \otimes \vc{v}^S- \mathbb{S}^S\right)\chi^S - \mathbb{S}^G\chi^G\\
			&+ \left( \varrho^s_1 \vc{v}^s_1 \otimes \vc{v}^s_1-\mathbb{S}^s_1 \right) \delta_{\Sigma_1}+ \left( \varrho^s_2 \vc{v}^s_2 \otimes \vc{v}^s_2-\mathbb{S}^s_2\right) \delta_{\Sigma_2} + \left( \varrho^s_3 \vc{v}^s_3 \otimes \vc{v}^s_3 -\mathbb{S}^s_3\right)  \delta_{\Sigma_3},
		\end{aligned}
	\end{equation}
	and the total energy flux as,
	\begin{align}
		\vc{J}_E 
		& :=\big( \varrho^L e^L \vc{v}^L - \mathbb{S}^L \vc{v}^L + \vc{J}_q^L\big) \chi^L + \big( \varrho^S e^S \vc{v}^S- \mathbb{S}^S \vc{v}^S + \vc{J}^S_q\big) \chi^S +\left( - \mathbb{S}^G \vc{v}^G + \vc{J}^G_q\right) \chi^G  \label{energy_flux}\\
		&  + \left( \varrho^s_1 e^s_1 \vc{v}^s_1 - \mathbb{S}^s_1 \vc{v}^s_1 + \vc{J}^s_{q,1}\right) \delta_{\Sigma_1}+\left(\varrho^s_2 e^s_2 \vc{v}^s_2 - \mathbb{S}^s_2 \vc{v}^s_2 + \vc{J}^s_{q,2}\right) \delta_{\Sigma_2} +\left(\varrho^s_3 e^s_3 \vc{v}^s_3 - \mathbb{S}^s_3 \vc{v}^s_3 + \vc{J}^s_{q,3}\right) \delta_{\Sigma_3}.\nonumber
	\end{align}
	The terms $\varrho^\xi \vc{v}^\xi \otimes \vc{v}^\xi,
	\; \varrho^s_i \vc{v}^s_i \otimes \vc{v}^s_i$ denote the transport of momentum in the bulk and at the interfaces respectively. Similarly, $\varrho^\xi e^\xi \vc{v}^\xi,  
	\; \varrho^s_i e^s_i \vc{v}^s_i$ denote the transport of energy.  The terms $\mathbb{S}^\xi,
	\mathbb{S}^s_i$ are the stress tensors that describe forces acting on a given element due to its surrounding. Therefore the terms
	$\mathbb{S}^\xi \vc{v}^\xi, 
	\;	\mathbb{S}^s_i \vc{v}^s_i$ yield the mechanical work made by these forces in any macroscopic element of the system. Finally, the terms $\vc{J}_q^\xi,
	\; \vc{J}_{q,i}^s$ denote the heat fluxes occurring in the system. In case of the bulk phases, 
	\begin{equation}
		\label{0.}
		\mathbb{S}^\xi = \sum_{k,l} s^\xi_{kl}(x) \; \vc{e}^k\otimes \vc{e}^l, \quad k,l=1,2,3, \quad \xi = L, G, S,
	\end{equation}
	where $\{\vc{e}^k\}$ is an orthonormal basis of $\mathbb{R}^3$. We do not distinguish the indices below or above since the scalar product of $\mathbb{R}^3$ is given by the metric $g_{kl} = \delta_{kl}$. In (\ref{0.}), $s^\xi_{kl}$ are functions $s^\xi_{kl}:D\to \mathbb{R}, D\subset \mathbb{R}^3$. On the other hand, $\mathbb{S}^s_i$ are the stress tensors in the tangent bundle of the surfaces $\Sigma_i$, i.e. 
	$$
	\forall \vc{x}\in\Sigma_i, \qquad \mathbb{S}^s_i(\vc{x}) \in T_{\vc{x}}\Sigma_i \otimes T_{\vc{x}}\Sigma_i, \quad i=1,2,3.
	$$
	Therefore, $ \mathbb{S}^s_i:\Sigma_i \to T\Sigma_i \otimes T\Sigma_i$.
	We also assume that all the stress tensors are symmetric (cf. \cite[Section 2.2.3]{shikh_bk}).
	The fluxes $\vc{J}_q^L, \vc{J}_q^G, \vc{J}_q^S \in\mathbb{R}^3$ are the standard vector fields where they are defined and $\vc{J}_{q,i}^s (\vc{x})\in T_{\vc{x}}\Sigma_i, i=1,2,3,$ for all $\vc{x}\in\Sigma_i$.
	
	Note that even if we neglect the contribution of the gas density and energy density in the conservation equations, there might be possible fluxes such as $ \vc{J}^G_q$.
	
	With these definitions, all the terms appearing in equations (\ref{momentum_flux}) and (\ref{energy_flux}) are meaningful in the sense of distributions. The conservation laws (\ref{mass_conv}), (\ref{momentum_conv}), (\ref{energy_conv}) are the equations describing the evolution of the full system, after suitable choices of the functions $\mathbb{S}, \mathbb{J}_q$ are made. Before describing these functions, we reformulate the above conservation laws  as a set of partial differential equations for the different magnitudes describing the evolution of the system (density, velocity, temperature), both in the bulk and at the interfaces. The conservation laws provide also a set of boundary conditions at the contact line $C$.
	\subsubsection{Conservation of mass}
	
	Equations \eqref{mass_conv}-\eqref{energy_conv} hold in the sense of distributions and imply corresponding laws in the bulk, on the interfaces and on the contact line. The derivation is rather standard, but in the following we present this derivation for \eqref{mass_conv} in detail for completeness and just state the corresponding results for the other conservation laws.

	Let us denote by $\Omega_L, \Omega_S$ the liquid and solid domains respectively, i.e. $\text{supp}\;\chi^L \equiv \Omega_L$ etc. In the  weak formulation, the mass conservation equation (\ref{mass_conv}) becomes, for $\varphi\in C^\infty_c( \mathbb{R}^3\times \mathbb{R})$,
	\begin{align}
		0 & = \int_0^\infty \int_{\mathbb{R}^3} {\Big[ \varrho^L \chi^L + \varrho^S \chi^S + \sum_{i=1}^3\varrho^s_i \delta_{\Sigma_i} \Big] \partial_t \varphi\; \mathrm{d}x \;\mathrm{d}t} \nonumber\\
		& \hspace{3cm} + \int_0^\infty \int_{\mathbb{R}^3} {\Big[ \varrho^L \vc{v}^L \chi^L +\varrho^S \vc{v}^S \chi^S +\sum_{i=1}^3 \varrho^s_i \vc{v}^s_i \delta_{\Sigma_i} \Big] \cdot \nabla_{\mathbf{x}}\varphi\; \mathrm{d}x \;\mathrm{d}t} \nonumber\\
		& = \int_0^\infty \mathrm{d}t {\Big[ \int_{\Omega_L}\varrho^L \;\partial_t \varphi\; \mathrm{d}x +\int_{\Omega_S}\varrho^S \;\partial_t \varphi\; \mathrm{d}x +\sum_{i=1}^3 \int_{\Sigma_i}\varrho^s_i \;\partial_t \varphi\; \mathrm{d}S \Big]} \nonumber\\
		& + \int_0^\infty \mathrm{d}t {\Big[ \int_{\Omega_L}\varrho^L \vc{v}^L \cdot \nabla_{\mathbf{x}}\varphi\; \mathrm{d}x +\int_{\Omega_S}\varrho^S \vc{v}^S \cdot \nabla_{\mathbf{x}}\varphi\; \mathrm{d}x +\sum_{i=1}^3 \int_{\Sigma_i}\varrho^s_i \vc{v}^s_i \cdot \nabla_{\mathbf{x}}\varphi\; \mathrm{d}S\Big] } =: I + II. \label{34}
	\end{align}
	We calculate $II$ using integration by parts with respect to the space variable. Recalling that $\varrho^S$ and $\varrho^L$ are constants, this gives 
	\begin{equation*}
		\begin{aligned}
			II 
			& = -\int_0^\infty\int_{\Omega_L}{\varrho^L\div\vc{v}^L \varphi \;\mathrm{d}x \;\mathrm{d}t} - \int_0^\infty\int_{\partial\Omega_L} {\varrho^L\vc{v}^L\cdot \vc{n} \;\varphi \;\mathrm{d}S\;\mathrm{d}t} -\int_0^\infty\int_{\Omega_S}{\varrho^S\div\vc{v}^S \varphi \;\mathrm{d}x \;\mathrm{d}t}\\
			& \quad \ + \int_0^\infty\int_{\partial\Omega_S} {\varrho^S\vc{v}^S\cdot \vc{n} \;\varphi \;\mathrm{d}S\;\mathrm{d}t} + \sum_{i=1}^3 \int_0^\infty\int_{\Sigma_i} {\underset{\left( \varrho^s_i \vc{v}^s_{i,\tau}\cdot \nabla_{\tau}\varphi + \varrho^s_i v^s_{i,n} \;\partial_n\varphi\right) }{\underbrace{\varrho^s_i \vc{v}^s_i \cdot \left( \nabla_{\tau}\varphi + \partial_n\varphi \;\vc{n}\right)}}\mathrm{d}S \;\mathrm{d}t}\\
			& = -\int_0^\infty\int_{\Omega_L}{\varrho^L\div \vc{v}^L \varphi \;\mathrm{d}x \;\mathrm{d}t} - \int_0^\infty\int_{\Sigma_1 \cup \Sigma_2} {\varrho^L v^L_n \;\varphi \;\mathrm{d}S\;\mathrm{d}t}+ \int_0^\infty\int_{\Sigma_2 \cup \Sigma_3} {\varrho^S v^S_n \;\varphi \;\mathrm{d}S\;\mathrm{d}t} \\
			& \quad \ -\sum_{i=1}^3\Big[ \int_0^\infty\int_{\Sigma_i}{\mathrm{div}_{\Sigma}\left( \varrho^s_i \vc{v}^s_{i,\tau}\right) \varphi\;\mathrm{d}S\;\mathrm{d}t}+ \int_0^\infty\int_{C}{\varrho^s_i \vc{v}^s_{i,\tau}\cdot \nu_i \;\varphi \;\mathrm{d} o\;\mathrm{d}t}\\
			& \qquad + \int_0^\infty\int_{\Sigma_i}{\varrho^s_i v^s_{i,n} \;\partial_n\varphi \;\mathrm{d}S\;\mathrm{d}t}\Big].
		\end{aligned}
	\end{equation*}
	Here we have denoted by $\vc{n}_i, i=1,2,$ the unit normal vectors to the surfaces $\Sigma_i, i=1,2$, inward with respect to the liquid domain and by $\vc{n}_3$ the unit normal vector to $\Sigma_3$ pointing into the gas. We have omitted here the subscript $i$ for the normal vectors for simplicity and will precise wherever it is required.
	Similarly, the unit vectors $\vc{\nu}_i$ are in the plane orthogonal to the contact line $C$ and also tangent to the interfaces $\Sigma_i, i=1,2,3$, pointing away from the contact line. The surface divergence $\mathrm{div}_\Sigma$ is defined in (\ref{surface_divergence.}).
	Also, we have used the specific form of the solid velocity field (\ref{solid_velocity}) in the above computation which yields $\div \vc{v}^S=0$.
	
	Furthermore, we use two transport equations (\ref{36a.}), (\ref{36b.}) in the following form, taking into account the choice of the normal vectors, which yield the boundary contributions at the interfaces $\Sigma_i(t)$ and at the moving contact line $C(t)$ respectively,
	\begin{equation}\label{36a}
		\int_0^\infty \int_{\Omega_L(t)}{\varrho^L \partial_t \varphi \;\mathrm{d}t\;\mathrm{d}x} = - 	\int_0^\infty \int_{\Omega_L(t)}{\partial_t\varrho^L \; \varphi \;\mathrm{d}t\;\mathrm{d}x} + \int_0^\infty \int_{\partial\Omega_L(t)}{\varrho^L v^s_{i,n} \;\varphi \;\mathrm{d}t\; \mathrm{d}S}, 
	\end{equation}
	similar for the solid domain $\Omega_S(t)$ and
	\begin{align}\label{36b}
		\int_0^\infty \int_{\Sigma_i(t)}{\varrho^s_i \partial_t^\Sigma \varphi \;\mathrm{d}t\;\mathrm{d}S}& = - \int_0^\infty \int_{\Sigma_i(t)}{\partial_t^\Sigma \varrho^s_i \varphi \;\mathrm{d}t\;\mathrm{d}S} - \int_0^\infty \int_{C}{\varrho^s_i \vc{v}^c\cdot \nu_i \;\varphi \;\mathrm{d}t \;\mathrm{d}o} \nonumber \\
		& \qquad - \int_0^\infty \int_{\Sigma_i(t)}{\varrho^s_i v^s_{i,n} \kappa_\Sigma\; \varphi\;\mathrm{d}t\;\mathrm{d}S}.
	\end{align}
	Here, $\kappa_\Sigma$ denotes the mean curvature of the interfaces and the normal derivative $\partial^\Sigma_t$ is defined in (\ref{normal_derivative.}).
	In addition, we assume that the normal component of the velocity of the interfaces in contact with the solid is equal to the normal velocity of the solid, i.e.
	\begin{equation}
		\label{normal_velocity}
		\left( \vc{v}^s_i - \vc{v}^S\right) \cdot \vc{n}_i =0 \quad \text{ at }\Sigma_i, \ \ i=2,3.
	\end{equation}
	This means that there is no flux of matter from the solid into the corresponding interfaces.
	
	Returning to the equation (\ref{34}) and integrating by parts with respect to the time variable and using \eqref{36a}-\eqref{normal_velocity}, we obtain,
	\begin{align}
		0 & = -\int\displaylimits_0^\infty\int\displaylimits_{\Omega_L}{\partial_t \varrho^L \varphi\; \mathrm{d}x\;\mathrm{d}t} +\int\displaylimits_0^\infty \int\displaylimits_{\partial\Omega_L}{\varrho^L v^s_{i,n}\varphi \;\mathrm{d}S\;\mathrm{d}t} -\int\displaylimits_0^\infty\int\displaylimits_{\Omega_S}{\partial_t \varrho^S \varphi\; \mathrm{d}x\;\mathrm{d}t} -\int\displaylimits_0^\infty \int\displaylimits_{\partial\Omega_S}{\varrho^S v^s_{i,n}\varphi \;\mathrm{d}S\;\mathrm{d}t} \nonumber\\
		& \quad + \sum_{i=1}^3 \int\displaylimits_0^\infty\int\displaylimits_{\Sigma_i}{\varrho^s_i \underset{\partial_t^\Sigma \varphi}{\underbrace{\left( \partial_t\varphi + v^s_{i,n}\;\partial_n\varphi\right)}}\; \mathrm{d}S\;\mathrm{d}t} -\int\displaylimits_0^\infty\int\displaylimits_{\Omega_L}{\varrho^L \div \vc{v}^L \varphi \;\mathrm{d}x \;\mathrm{d}t}  - \int\displaylimits_0^\infty\int\displaylimits_{\Sigma_1 \cup \Sigma_2} {\varrho^L v^L_n \;\varphi \;\mathrm{d}S\;\mathrm{d}t} \nonumber\\
		& \quad  + \int\displaylimits_0^\infty\int\displaylimits_{\Sigma_2 \cup \Sigma_3} {\varrho^S v^S_n \;\varphi \;\mathrm{d}S\;\mathrm{d}t} - \sum_{i=1}^3\Big[ \int\displaylimits_0^\infty\int\displaylimits_{\Sigma_i}{\mathrm{div}_{\Sigma}\left( \varrho^s_i \vc{v}^s_{i,\tau}\right) \varphi\;\mathrm{d}S\;\mathrm{d}t} + \int\displaylimits_0^\infty\int\displaylimits_{C}{\varrho^s_i \vc{v}^s_{i,\tau}\cdot \nu_i \;\varphi \;\mathrm{d} o\;\mathrm{d}t} \Big]\nonumber\\
		& = -\int\displaylimits_0^\infty\int\displaylimits_{\Omega_L}{ \varrho^L \div \vc{v}^L \varphi\; \mathrm{d}x\;\mathrm{d}t}- \int\displaylimits_0^\infty\int\displaylimits_{\Sigma_1 \cup \Sigma_2} {\varrho^L\left( v^L_n - v^s_{i,n}\right)  \varphi \;\mathrm{d}S\mathrm{d}t} + \int\displaylimits_0^\infty\int\displaylimits_{\Sigma_2 \cup \Sigma_3} {\varrho^S\left( v^S_n - v^s_{i,n}\right)  \varphi \;\mathrm{d}S\mathrm{d}t} \nonumber\\
		& \quad -\sum_{i=1}^3\Big[ \int\displaylimits_0^\infty\int\displaylimits_{\Sigma_i}{\left( \partial_t^\Sigma\varrho^s_i + \kappa_\Sigma \varrho^s_i v^s_{i,n}\right) \varphi\; \mathrm{d}S\;\mathrm{d}t}- \int\displaylimits_0^\infty\int\displaylimits_C {\varrho^s_i \vc{v}^c \cdot \nu_i \;\varphi \;\mathrm{d}o \;\mathrm{d}t} \nonumber\\
		& \quad \quad - \int\displaylimits_0^\infty\int\displaylimits_{\Sigma_i}{\mathrm{div}_{\Sigma}\left( \varrho^s_i \vc{v}^s_{i,\tau}\right) \varphi\;\mathrm{d}S\;\mathrm{d}t}+ \int\displaylimits_0^\infty\int\displaylimits_{C}{\varrho^s_i \vc{v}^s_{i,\tau}\cdot \nu_i \;\varphi \;\mathrm{d} o\;\mathrm{d}t} \Big]. \label{63}
	\end{align}
	Therefore, mass conservation in the liquid and at the interfaces $\Sigma_i$ become,
	\begin{equation}
		\label{incom}
		\div \vc{v}^L=0 \quad \text{ in } \Omega_L,
	\end{equation}
	\begin{equation}
		\label{37}
		\begin{aligned}
			&	\partial_t^{\Sigma} \varrho^s_i + \mathrm{div}_{\Sigma} \left( \varrho^s_i \vc{v}^s_{i,\tau}\right) + \kappa_\Sigma \varrho^s_i v^s_{i,n} + \varrho^L (\delta_{1,i}+\delta_{2,i}) \left( v^L_{n} - v^s_{i,n}\right)  =0 \quad \text{at } \Sigma_i, i=1,2,3.
		\end{aligned}
	\end{equation}
	Note that the mass balance equation gives no condition for the gas in the limit $\varrho^G=0$. Also there is no contribution of the solid density $\varrho^S$ to the interface equation (\ref{37}) due to the condition (\ref{normal_velocity}).
	
	Equation (\ref{37}) for $i=2$, that is the mass condition at the liquid-solid interface, is a generalization of the standard boundary condition $v^L_n =v^s_{2,n}=v^S_{n} $ for an immiscible fluid.
	
	We also obtain the boundary condition for the function $\varrho^s_i$ at the contact line $C$, derived from the mass balance equation (\ref{mass_conv}), or equivalently from (\ref{63}) as,
	\begin{equation}
		\label{mass_cont}
		\sum_{i=1}^3 \varrho^s_i \left( \vc{v}^s_{i,\tau} - \vc{v}^c\right)\cdot \vc{\nu}_i =0 \quad \text{ on } C.
	\end{equation}

	\subsubsection{Conservation of momentum}
	Similarly, we obtain from the momentum balance equation (\ref{momentum_conv}) the following
	equations in the bulk and at the interfaces $\Sigma_i$,
	\begin{equation}
		\label{NS}
		\partial_t \big( \varrho^\xi\vc{v}^\xi \big) + \div \big( \varrho^\xi \vc{v}^\xi \otimes \vc{v}^\xi - \mathbb{S}^\xi\big) = 0\, \quad \mbox{ in } \Omega_\xi,\ \xi=L,S,
	\end{equation}
	
	\begin{equation}
		\label{pressure_gas}
		\div \mathbb{S}^G=0\qquad \mbox{ in } \Omega_G,
	\end{equation}
	\begin{equation}
		\begin{aligned}
			\label{17...}
			&	\partial^\Sigma_t (\varrho^s_i \vc{v}^s_i) + \mathrm{div}_\Sigma (\varrho^s_i \vc{v}^s_i \otimes \vc{v}^s_{i,\tau} - \vm{S}^s_i) + \kappa_\Sigma \varrho^s_i \vc{v}^s_i v^s_{i,n}\\
			& \hspace{1.6cm} + \varrho^L \vc{v}^L\left( v^L_n- v^s_{i,n}\right) =\big(\vm{S}^L - \vm{S}^\xi\big)\cdot{\vc{n}_i} \quad \text{ on } \Sigma_i, \ (i,\xi) \in \{(1,G),(2,S)\},
		\end{aligned}
	\end{equation}
	
	\begin{equation}
		\label{17..}
		\partial^\Sigma_t (\varrho^s_3 \vc{v}^s_3) + \mathrm{div}_\Sigma (\varrho^s_3 \vc{v}^s_3 \otimes \vc{v}^s_{3,\tau} - \vm{S}^s_3) + \kappa_\Sigma \varrho^s_3 \vc{v}^s_3 v^s_{3,n}  =\big(\vm{S}^G - \vm{S}^S\big)\cdot{\vc{n}_3} \quad \text{ at } \Sigma_3.
	\end{equation}
	Note that the contribution of the flux of momentum from the solid does not appear in equation (\ref{17...}), $(i, \xi)=(2,S)$ and (\ref{17..}) due to the continuity assumption of the normal velocity (\ref{normal_velocity}).
	
	The surface momentum equation (\ref{17...}) with $(i, \xi) = (2,S)$ is analogous to the standard capillary equation
	\begin{equation*}
		\mathrm{div}_\Sigma \vm{S}^s+\big(\vm{S}^L - \vm{S}^S\big)\cdot \vc{n} =0,
	\end{equation*}
	which can be obtained from (\ref{17...}) by taking $\varrho^s_2=0$ and $v^L_n = v^s_{2,n}$.
	
	Further, the momentum balance equation (\ref{momentum_conv}) implies at the contact line $C$,
	\begin{equation}
		\label{64}
		\sum_{i=1}^3 \Big[ \varrho^s_i \vc{v}^s_i \otimes\left( \vc{v}^s_{i,\tau} - \vc{v}^c\right) - \mathbb{S}^s_i\Big] \cdot \vc{\nu}_i =0 \qquad \mbox{ on } C.
	\end{equation}
	Condition (\ref{64}) is equivalent to Young's law in the standard case.
	
	\subsubsection{Conservation of energy}
	
	Next, the energy equations in the bulk and at the interfaces follow from the energy conservation law (\ref{energy_conv}),
	\begin{equation}
		\label{24}
		\partial_t (\varrho^\xi e^\xi) + \div (\varrho^\xi e^\xi \vc{v}^\xi -\mathbb{S}^\xi \cdot\vc{v}^\xi + \vc{J}^\xi_q)=0 \quad \mbox{ in } \Omega_\xi, \ \xi=L,S,
	\end{equation}
	
	\begin{equation}
		\label{22.}
		\div \vc{J}^G_q=\div \left( \mathbb{S}^G \cdot\vc{v}^G\right) \quad \mbox{ in } \Omega_G,
	\end{equation}
	\begin{equation}
		\label{24.1}
		\begin{aligned}
			&\partial^\Sigma_t (\varrho^s_i e^s_i) + \mathrm{div}_\Sigma \left( \varrho^s_i e^s_i \vc{v}^s_{i,\tau} - \vm{S}^s_i \cdot \vc{v}^s_i + \vc{J}^s_{i,q} \right) + \kappa_\Sigma \varrho^s_i e^s_i v^s_{i,n} + \varrho^L e^L \left( v^L_n - v^s_{i,n}\right) \\
			&\hspace{1.5cm}= \left(\vm{S}^L\cdot \vc{v}^L - \vc{J}_q^L -\vm{S}^\xi\cdot \vc{v}^\xi + \vc{J}_q^\xi\right)\cdot{\vc{n}_i}\qquad \text{ on } \ \Sigma_i, \ (i,\xi) \in \{(1,G),(2,S)\},
		\end{aligned}
	\end{equation}
	
	\begin{equation}
		\label{24.3}
		\begin{aligned}
			&\partial^\Sigma_t (\varrho^s_3 e^s_3) + \mathrm{div}_\Sigma \left( \varrho^s_3 e^s_3 \vc{v}^s_{3,\tau} - \vm{S}^s_3 \cdot \vc{v}^s_3 + \vc{J}^s_{3,q} \right) + \kappa_\Sigma \varrho^s_3 e^s_3 v^s_{3,n} \\
			& \hspace{4.5cm}= \left(\vm{S}^G \cdot \vc{v}^G - \vc{J}_q^G -\vm{S}^S \cdot \vc{v}^S + \vc{J}_q^S\right)\cdot{\vc{n}_3} \qquad \text{ at } \Sigma_3.
		\end{aligned}
	\end{equation}
	Recall that $e^\xi, e^s_i$ in the above equations are the total energies.
	
	From (\ref{24}), with the help of the balance equation for the kinetic energy in the liquid, which can be obtained by multiplying the momentum equation (\ref{NS}) by $\vc{v}^L$, together with the mass equation (\ref{incom}) as (similarly for the solid),
	\begin{equation*}
		\partial_t \Big(\varrho^L\frac{|\vc{v}^L|^2}{2} \Big) + \div \left( \varrho^L \frac{|\vc{v}^L|^2}{2}\vc{v}^L \right) = \vc{v}^L\cdot\div \mathbb{S}^L,
	\end{equation*}
	one obtains the balance equations for the internal energy $u^L$ and $u^S$ as,
	\begin{equation}
		\label{22}
		\partial_t (\varrho^\xi u^\xi) + \div \left( \varrho^\xi u^\xi \vc{v}^\xi + \vc{J}^\xi_q \right) = \vm{S}^\xi : \nabla_{\vc{x}} \vc{v}^\xi\quad \mbox{ in } \Omega_\xi, \ \xi=L,S.
	\end{equation}	
	Similarly, subtracting the surface mass equation (\ref{37}) multiplied by $\frac{1}{2}|\vc{v}^s_i|^2 $ from the surface momentum equation (\ref{17..}) multiplied by $\vc{v}^s_i$ yields the balance equations for the kinetic energy on the interfaces, 
	\begin{equation*}
		\begin{aligned}
			&\partial_{t}^\Sigma\Big( \varrho^s_i \frac{|\vc{v}^s_i|^2}{2}\Big) + \mathrm{div}_\Sigma \Big( \varrho^s_i \vc{v}^s_{i,\tau} \frac{|\vc{v}^s_i|^2}{2}\Big) + \kappa_\Sigma \varrho^s_i v^s_{i,n} \frac{|\vc{v}^s_i|^2}{2} + \varrho^L \vc{v}^s_i \cdot \Big( \vc{v}^L - \frac{1}{2}\vc{v}^s_i\Big) \left(v^L_n - v^s_{i,n}\right) \\
			& \hspace{3.5cm} - \vc{n}_i \cdot \big( \mathbb{S}^L - \mathbb{S}^\xi\big) \cdot \vc{v}^s_i=\vc{v}^s_i \cdot \mathrm{div}_\Sigma \; \mathbb{S}^s_i \qquad \text{on } \Sigma_i, \ (i, \xi)\in \{(1,G),(2,S)\};
		\end{aligned}
	\end{equation*}
	Indeed, with the standard formula $\nabla \cdot \left( a\otimes b\right) = (\nabla \cdot a) b + (a\cdot \nabla )b$, we can compute the convective term as,
	\begin{equation*}
		\begin{aligned}
			& \ -\frac{1}{2}|\vc{v}^s_i|^2 \mathrm{div}_\Sigma\left( \varrho^s_i \vc{v}^s_{i,\tau}\right) + \vc{v}^s_i \cdot \mathrm{div}_\Sigma \left( \varrho^s_i \vc{v}^s_i \otimes \vc{v}^s_{i,\tau}\right) \\
			&= -\frac{1}{2}|\vc{v}^s_i|^2 \mathrm{div}_\Sigma\left( \varrho^s_i \vc{v}^s_{i,\tau}\right) + \vc{v}^s_i \cdot\left[ \mathrm{div}_\Sigma \left( \varrho^s_i \vc{v}^s_i\right) \vc{v}^s_{i,\tau} +  \varrho^s_i \left( \vc{v}^s_i\cdot \nabla_\tau \right) \vc{v}^s_{i,\tau}\right] \\
			&= \frac{1}{2}|\vc{v}^s_i|^2 \mathrm{div}_\Sigma\left( \varrho^s_i \vc{v}^s_{i,\tau}\right) + \varrho^s_i\vc{v}^s_i \cdot \left( \vc{v}^s_i\cdot \nabla_\tau \right) \vc{v}^s_{i,\tau} = \mathrm{div}_\Sigma \Big( \varrho^s_i \vc{v}^s_{i,\tau} \frac{|\vc{v}^s_i|^2}{2}\Big),
		\end{aligned}
	\end{equation*}
	which is appearing in the above equation for the surface kinetic energy. We also find the following equation on the liquid-gas interface,
	\begin{equation*}
		\begin{aligned}
			&\partial_{t}^\Sigma\Big( \varrho^s_3 \frac{|\vc{v}^s_3|^2}{2}\Big) + \mathrm{div}_\Sigma \Big( \varrho^s \vc{v}^s_{3, \tau} \frac{|\vc{v}^s_3|^2}{2}\Big)+ \kappa_\Sigma\varrho^s_3 v^s_{3,n} \frac{|\vc{v}^s_3|^2}{2}
			\\
			& \hspace{4.9cm} - \vc{n}_3 \cdot \left( \mathbb{S}^G - \mathbb{S}^S\right) \cdot \vc{v}^s_3 =  \vc{v}^s_3\cdot \mathrm{div}_\Sigma \; \mathbb{S}^s_3 \qquad \text{ on } \Sigma_3.
		\end{aligned}
	\end{equation*}
	
	Then the following balance equations for the internal energy $u^s_i$ are obtained from (\ref{24.1})-(\ref{24.3}), with the help of the above kinetic energy equations,\\
	\begin{equation}
		\label{11.}
		\begin{aligned}
			& \partial^\Sigma_t (\varrho^s_i u^s_i) + \mathrm{div}_\Sigma \left( \varrho^s_i u^s_i \vc{v}^s_{i,\tau} + \vc{J}^s_{q,i} \right) + \kappa_\Sigma\varrho^s_i u^s_i v^s_{i,n} + \varrho^L u^L \left( v^L_n - v^s_{i,n}\right) = \vm{S}^s_i : \nabla_{\tau}\vc{v}^s_i \\
			& +\vc{n}_i\cdot \mathbb{S}^L\cdot(\vc{v}^L - \vc{v}^s_i) - \vc{n}_i \cdot \mathbb{S}^\xi\cdot(\vc{v}^\xi - \vc{v}^s_i)- \left( \vc{J}_q^L - \vc{J}_q^\xi\right)\cdot{\vc{n}_i}\\
			& - \frac{1}{2}\varrho^L (v^L_n - v^s_{i,n}) (\vc{v}^L-\vc{v}^s_i)^2 \qquad \text{ on } \Sigma_i, \ (i, \xi)\in \{(1,G),(2,S)\},
		\end{aligned}
	\end{equation}
	and 
	\begin{equation}
		\label{GS_energy}
		\begin{aligned}
			& \partial^\Sigma_t (\varrho^s_3 u^s_3) + \mathrm{div}_\Sigma \left( \varrho^s_3 u^s_3 \vc{v}^s_{3,\tau} + \vc{J}^s_{q,3} \right) + \kappa_\Sigma\varrho^s_3 u^s_3 v^s_{3,n} = \vm{S}^s_3 : \nabla_{\tau}\vc{v}^s_3  \\
			& \hspace{1.2cm}+ \vc{n}_3\cdot \mathbb{S}^G\cdot(\vc{v}^G - \vc{v}^s_3) - \vc{n}_3 \cdot \mathbb{S}^S\cdot(\vc{v}^S - \vc{v}^s_3)- \left( \vc{J}_q^G - \vc{J}_q^S\right)\cdot{\vc{n}_3} \qquad \text{ on } \Sigma_3.
		\end{aligned}
	\end{equation}
	The last term of equation (\ref{11.}) appears due to the difference of the velocities of the interface and the liquid which transforms into the internal energy.
	
	The energy equation at the contact line $C$ becomes, from the conservation equation (\ref{energy_conv}),
	\begin{equation}
		\label{27}
		\sum_{i=1}^3 \varrho^s_i e^s_i \left( \vc{v}^s_{i,\tau} - \vc{v}^c\right) \cdot \vc{\nu}_i + \left( - \mathbb{S}^s_i \cdot \vc{v}^s_i + \vc{J}^s_{q,i}\right)\cdot \vc{\nu}_i =0 \qquad \text{ on } C.
	\end{equation}
	The above equation (\ref{27}) can be rewritten equivalently, with the help of the mass and momentum equations at the contact line (\ref{mass_cont}) and (\ref{64}) as,
	\begin{equation}
		\label{energy_cont}
		\sum_{i=1}^3\Big(\varrho^s_i \Big( u^s_i + \frac{|\vc{v}^s_i - \vc{v}^c|^2}{2}\Big) (\vc{v}^s_{i,\tau} - \vc{v}^c)- \mathbb{S}^s_i \cdot(\vc{v}^s_{i,\tau} - \vc{v}^c) + \vc{J}_{q,i}^s\Big)\cdot \vc{\nu}_i =0 \qquad \text{ on } C.
	\end{equation}
	
	\subsection{Constitutive equations}
	
	We now use Onsager's method to derive relations between the generalized fluxes and generalized forces. To this end, we first write the equation for the balance of entropy.
	\subsubsection{Entropy balance}
	The general form of the equation satisfied by the specific entropy is,
	
	\begin{equation}
		\label{entropy}
		\partial_t \tilde{s} + \div \left( \tilde{s} \vc{v} + \vc{J}_s\right) = \sigma,
	\end{equation}
	where the total specific entropy is expressed as, given the previous expression for the energy density (i.e. neglecting the contribution of the gas),
	$$
	\tilde{s}:= \varrho^L s^L \chi^L + \varrho^S s^S \chi^S + \varrho^s_1 s^s_1 \delta_{\Sigma_1}+ \varrho^s_2 s^s_2 \delta_{\Sigma_2} + \varrho^s_3 s^s_3 \delta_{\Sigma_3},
	$$
	the total entropy flux as,
	$$\vc{J}_s :=\vc{J}_s^L \chi^L + \vc{J}^S_s \chi^S +\vc{J}^G_s \chi^G + \vc{J}^s_{s,1} \delta_{\Sigma_1}+\vc{J}^s_{s,2} \delta_{\Sigma_2} + \vc{J}^s_{s,3} \delta_{\Sigma_3},
	$$
	and the total entropy production as,
	$$\sigma :=\sigma^L \chi^L +\sigma^S \chi^S +\sigma^G \chi^G + \sigma^s_1 \delta_{\Sigma_1}+\sigma^s_2 \delta_{\Sigma_2} + \sigma^s_3 \delta_{\Sigma_3}.
	$$
	Here $s^\xi, s^s_i$ are the specific entropies per unit mass, $\vc{J}^\xi, \vc{J}^s_i$ are the entropy fluxes per unit mass and $\sigma^\xi, \sigma^s_i$ entropy production per unit volume in the bulk and at the interfaces respectively.
	
	Note that, even if we neglect the contribution of the gas density and energy density (hence, entropy density) in the conservation equations, there might be possible fluxes such as $\vc{J}^G_q, \vc{J}^G_s$. The entropy fluxes are functions of the energy fluxes, heat fluxes etc, discussed previously.
	
	The above balance equation (\ref{entropy}) yields the classical entropy equations in the bulk, that is,
	\begin{equation}
		\label{21}
		\partial_t (\varrho^\xi s^\xi) + \div \big( \varrho^\xi s^\xi \vc{v}^\xi + \vc{J}_s^\xi \big) = \sigma^\xi \qquad \mbox{ in } \Omega_{\xi}, \ \xi=L,S,
	\end{equation}
	
	\begin{equation}
		\label{21a}
		\div \; \vc{J}_s^G = \sigma^G \qquad \mbox{ in } \Omega_G,
	\end{equation}
	and at the interfaces, 
	\begin{equation}
		\label{21b}
		\begin{aligned}
			\partial^\Sigma_t (\varrho^s_i s^s_i) & + \mathrm{div}_\Sigma\left( \varrho^s_i s^s_i \vc{v}^s_{i,\tau} + \vc{J}^s_{s,i} \right) + \kappa_\Sigma \varrho^s_i s^s_i v^s_{i,n}+ \varrho^L s^L\left( v^L_n - v^s_{i,n}\right)\\
			& \hspace{2cm} + \big( \vc{J}_s^L - \vc{J}_s^\xi\big)\cdot{\vc{n}_i}  = \sigma^s_i  \qquad \mbox{ on } \Sigma_i, \; (i,\xi)\in \{(1,G),(2,S)\}, 
		\end{aligned}
	\end{equation}
	and
	\begin{equation}\label{21c}
		\partial^\Sigma_t (\varrho^s_3 s^s_3) + \mathrm{div}_\Sigma\left( \varrho^s_3 s^s_3 \vc{v}^s_{3,\tau} + \vc{J}^s_{s,3} \right) + \kappa_\Sigma \varrho^s_1 s^s_3 v^s_{3,n}+ \left( \vc{J}_s^G - \vc{J}_s^S\right)\cdot{\vc{n}_3}= \sigma^s_3  \qquad \mbox{ on } \Sigma_3.
	\end{equation}
	Finally, at the contact line we obtain,
	\begin{equation}
		\label{CL_entropy}
		\sum_{i=1}^3\Big(\varrho^s_i s^s_i (\vc{v}^s_{i,\tau} - \vc{v}^c)+ \vc{J}^s_{s,i} \Big) \cdot \vc{\nu}_i= 0 \qquad \mbox{ on } C.
	\end{equation}
	
	The goal here is to obtain the entropy productions $\sigma^\xi, \sigma^s_i$ from the energy equations and the Gibbs relation (using the hypothesis of local equilibrium).
	
	Classical thermodynamics proposes the relation 
	$dS = \frac{dQ}{T}$ between entropy and heat flux.
	Combining this with the assumption of local equilibrium, it is natural to assume that 
	\begin{equation}
		\label{fluxes}
		\vc{J}_s^\xi = \frac{\vc{J}_q^\chi}{T}, \qquad \vc{J}_{s,i}^s = \frac{\vc{J}_{q,i}^s}{T}, \qquad \xi = L, S, G, \ \ i=1, 2, 3.
	\end{equation}
	In order to compute the production of entropy $\sigma^L, \sigma^S, \sigma^G, \sigma^s_i$, we use the the general Gibbs equation (cf. \cite{DGM}),
	\begin{equation*}
		\mathrm{d}(\varrho s) = \frac{1}{T}\left( \mathrm{d}(\varrho u) - \mu \;\mathrm{d}\varrho\right)
	\end{equation*}
	where $T$ is the absolute temperature and $\mu$ is the chemical potential per unit mass and the thermodynamic pressure $p$ is related to the chemical potential by
	\begin{equation*}
		\label{}
		\varrho s = \frac{1}{T} \left( \varrho u - \varrho\mu + p\right).
	\end{equation*}
	Taking into account	the assumptions that the densities of the liquid and the solid are constants,
	the above Gibbs relation gives,
	\begin{align}
		\mathrm{d}s^\xi&= \frac{1}{T}\mathrm{d}u^\xi, \quad \xi=L,S,\label{gibbs} \\
		\mathrm{d}(\varrho^s_i s^s_i) &= \frac{1}{T}\left( \mathrm{d}(\varrho^s_i u^s_i) - \mu^s_i \;\mathrm{d}\varrho^s_i\right) , \quad i=1,2,3,\label{gibbs1}
	\end{align}
	together with
	\begin{equation}
		\label{20.}
		\varrho^\xi s^\xi = \frac{1}{T} \left( \varrho^\xi u^\xi - \varrho^\xi \mu^\xi + p^\xi\right), \quad \varrho^s_i s^s_i = \frac{1}{T} \left( \varrho^s_i u^s_i - \varrho^s_i \mu^s_i + p^s_i\right), \quad \xi=L,S, \ i=1,2,3.
	\end{equation} 
	
	The meaning of the relations (\ref{gibbs}), (\ref{gibbs1}) are precisely that the entropy density of bulk and the interfaces are functions of (internal) energy density and mass density only. Thus we have the constitutive relations
	\begin{equation}
		\label{gibbs2}
		\varrho^\xi s^\xi = h^\xi(\varrho^\xi u^\xi), \qquad \varrho^s_i s^s_i = h^s_i(\varrho^s_i u^s_i, \varrho^s_i), \qquad \xi = L,S,  \quad i=1,2,3,
	\end{equation}
	for some functions $h^{\xi}$, $h_i^s$	and the absolute temperature and the chemical potential in the respective phases are defined as,
	\begin{equation}
		\label{gibbs3}
		\frac{1}{T} = \frac{\partial h^\xi}{\partial (\varrho^\xi u^\xi)} = \frac{\partial h^s_i}{\partial (\varrho^s_i u^s_i)}, 
		\qquad -\frac{\mu^s_i}{T} = \frac{\partial h^s_i}{\partial \varrho^s_i}, \qquad \xi = L,S,  \quad i=1,2,3.
	\end{equation}
	
	Here the strategy is to obtain expressions for the entropy productions $\sigma$ from the entropy equations (\ref{21})-(\ref{21c})
	by using the Gibbs relations $\varrho s = h(\varrho, \varrho u)$ in each phase.
	Writing the convective derivatives $\left( \partial_{t} + \vc{v}\cdot \nabla \right)(\cdot) $ with the help of the relation (\ref{gibbs3}) and then plugging in the mass and energy equations for $ \varrho, \varrho u$ yields $\sigma$ in terms of the generalized flux-force pairs, where one further employs Onsager principle.
	
	In the bulk, we therefore write down the entropy production (\ref{21}) with the help of Gibbs relation (\ref{gibbs}) and (\ref{20.}) (or, (\ref{gibbs2})-(\ref{gibbs3})) and the flux relation (\ref{fluxes}), using energy equation (\ref{22}) and the incompressibility condition (\ref{incom}), as
	\begin{align}
		\sigma^L &= \partial_t \big( \varrho^L s^L\big) + \div \big( \varrho^L s^L\vc{v}^L \big) + \div \vc{J}^L_s \nonumber\\
		&=\frac{1}{T}\big( \partial_t \left( \varrho^L u^L\right) +\vc{v}^L \cdot \nabla_{\mathbf{x}}\left( \varrho^L u^L\right) \big) +\varrho^L s^L \div \vc{v}^L + \div \Big( \frac{\vc{J}^L_q}{T} \Big) \nonumber\\
		& = \frac{1}{T}\Big( \mathbb{S}^L : \nabla_{\mathbf{x}}\;\vc{v}^L - \div \vc{J}^L_q \Big) - \frac{1}{T}\varrho^L u^L \div \vc{v}^L+\varrho^L s^L \div \vc{v}^L+ \div \Big( \frac{\vc{J}^L_q}{T} \Big) \nonumber\\
		&= \vc{J}^L_q \cdot \nabla_{\mathbf{x}}\Big( \frac{1}{T}\Big) +  \widetilde{\mathbb{S}^L}:\frac{1}{T}\mathbb{D}\vc{v}^L \qquad \text{ in } \Omega_L, \label{25} \\
		\sigma^S &= \vc{J}^S_q \cdot \nabla_{\mathbf{x}}\Big( \frac{1}{T}\Big) \qquad \text{ in } \Omega_S, \label{25.1} \\
		\sigma^G& = \vc{J}^G_q \cdot \nabla_{\mathbf{x}}\Big( \frac{1}{T}\Big)+  \widetilde{\mathbb{S}^G}:\frac{1}{T}\mathbb{D}\vc{v}^G + \left( p^G - P^G\right) \frac{1}{T}\div\vc{v}^G \qquad \text{ in } \Omega_G. 	\label{25.2}
	\end{align}
	We denote here and below by $ \widetilde{\mathbb{S}^\xi} := (\mathbb{S}^\xi + P^\xi \mathbb{I}), \ \widetilde{\mathbb{S}^s_i} := (\mathbb{S}^s_i + P^s_i \left( \mathbb{I} - \vc{n}_i\otimes \vc{n}_i\right)) $ the traceless part of the stress tensors where $P^\xi := -\frac{1}{3} \mathrm{tr}\;\mathbb{S}^\xi, \ P^s_i := -\frac{1}{2} \mathrm{tr}\;\mathbb{S}^s_i$ are the mechanical pressures in the bulk and at the interfaces respectively. The factor $\frac{1}{3}$ in the mechanical pressure is appearing because of the three-dimensional bulk phases and the factor $\frac{1}{2}$ appears due to the two-dimensional interfaces. We also denote by $\mathbb{D}\vc{v}^\xi := \frac{1}{2}\left( \nabla_{\mathbf{x}} + \nabla_{\mathbf{x}}^T\right)\vc{v}^\xi, \ \mathbb{D}\vc{v}^s_i := \frac{1}{2}\left( \mathbb{I} - \vc{n}_i\otimes \vc{n}_i\right) \left( \nabla_{\tau} + \nabla_{\tau}^T\right)\vc{v}^s_i \left( \mathbb{I} - \vc{n}_i\otimes \vc{n}_i\right)$ the symmetric gradients. 
	In the deduction of (\ref{25.1}), we have used the fact that $\mathbb{D}\vc{v}^S =0$ due to the rigid body motion (\ref{solid_velocity}).
	
	Similarly we obtain from (\ref{21b})-(\ref{21c}) together with (\ref{gibbs1})-(\ref{20.}), (\ref{fluxes}) combining (\ref{11.}) and (\ref{37}),
	at $\Sigma_i, \ (i,\xi) \in \{(1,G), (2,S)\}$,
	\begin{align}
		\sigma^s_i = & \ \vc{J}^s_{q,i}\cdot \nabla_\tau\Big( \frac{1}{T}\Big)+ \frac{1}{T}\widetilde{\mathbb{S}^s_i}: \mathbb{D}\vc{v}^s_{i,\tau} +\frac{1}{T}\left( p^s_i -P^s_i\right) \left( \mathrm{div}_\Sigma \vc{v}^s_{i,\tau} +  \kappa_\Sigma v^s_{i,n}\right)  \nonumber \\
		& + \frac{1}{T}\varrho^L(v^L_n - v^s_{i,n}) \Big[ \left( \mu^s_i - \mu^L\right) + \frac{1}{\varrho^L }\;{\vc{n}}_i\cdot \left( \mathbb{S}^L+p^L \mathbb{I}\right) \cdot \vc{n}_i\Big] \label{25.}\\
		& + \frac{1}{2T} \Big[ \Big( \mathbb{S}^L + \mathbb{S}^\xi\Big)\cdot \vc{n}_i\Big] _{\tau} \cdot \big( \vc{v}^L - \vc{v}^\xi\big) _\tau - \frac{1}{T} \left[ \left( \mathbb{S}^L - \mathbb{S}^\xi\right)\cdot \vc{n}_i\right] _{\tau}\cdot \big( \vc{v}^s_i - \frac{1}{2}\big( \vc{v}^L+\vc{v}^\xi\big) \big) _\tau, \nonumber
	\end{align}
	and at $\Sigma_3$,
	\begin{align}
		\sigma^s_3 & = \vc{J}^s_{q,3}\cdot \nabla_\tau\Big( \frac{1}{T}\Big)+ \frac{1}{T}\widetilde{\mathbb{S}^s_3}: \mathbb{D}\vc{v}^s_{3, \tau} + \frac{1}{T}\left( p^s_3 -P^s_3\right) \left( \mathrm{div}_\Sigma \vc{v}^s_{3,\tau}+ \kappa_\Sigma v^s_{3,n} \right) \label{25...}\\
		& \ + \frac{1}{2T} \left[ \left( \mathbb{S}^G + \mathbb{S}^S\right)\cdot \vc{n}_3\right] _{\tau} \cdot \left( \vc{v}^G - \vc{v}^S\right) _\tau - \frac{1}{T} \big[ \left( \mathbb{S}^G - \mathbb{S}^S\right)\cdot \vc{n}_3\big] _{\tau}\cdot \big( \vc{v}^s_3 - \frac{1}{2}\big( \vc{v}^G+\vc{v}^S\big) \big) _\tau. \nonumber
	\end{align}
	
	In principle, the mechanical pressures $P^\xi, P^s_i$ could be different from the thermodynamical pressures $p^\xi, p^s_i$. In the case of the liquid, this difference is irrelevant as we assume that it is incompressible. Therefore $P^L$ is just Lagrange multiplier that guarantees that $\div \vc{v}^L=0$ is satisfied. On the other hand, the density of the interfaces could change and therefore the differences $(P^s_i-p^s_i)$ might be relevant.
	
	\subsubsection{Onsager's principle }
	
	At this point, we will use Onsager's idea to relate the generalized forces and fluxes which says that,
	\begin{itemize}
		\item 
		There is a linear relation between generalized forces and generalized fluxes. 
		\item
		The generation of entropy must be a non-negative  quadratic form.
		\item The matrix of the coefficients is symmetric in the absence of the magnetic field.
	\end{itemize}
	
	With the help of this principle, one obtains different constitutive equations from the entropy production rate. For example, on the right hand side of (\ref{25}), $\vc{J}^L_q$ represents the heat flux and $\nabla_{\mathbf{x}}\left( \frac{1}{T}\right) $ is the force, generated by the gradient of the temperature. Similarly, for the second term, $\widetilde{\mathbb{S}^L}$ denotes the viscous stress flux and $\mathbb{D}\vc{v}^L $ is the force due to the velocity gradient. Hence, Onsager's principle yields in $\Omega_L$,
	\begin{equation*}
		\begin{aligned}
			\vc{J}^L_q &= a^L \;\nabla_{\mathbf{x}}\Big( \frac{1}{T}\Big) + a_{T,v} \;\frac{1}{T}\mathbb{D}\vc{v}^L,\qquad 
			\widetilde{\mathbb{S}^L} &= a_{v,T}\; \nabla_{\mathbf{x}}\Big( \frac{1}{T}\Big) + \frac{b^L}{T} \mathbb{D}\vc{v}^L.
		\end{aligned}	
	\end{equation*}
	Here we further employ the assumption that the effects of the cross-diffusion terms are negligible, i.e. $a_{T,v} = a_{v,T} =0$.
	In the present scenario, we therefore obtain, from (\ref{25}), (\ref{25.2}), the usual form of the stress tensors in the bulk, 
	\begin{equation}
		\label{2..}
		\mathbb{S}^L = - P^L \mathbb{I} + \frac{b^L}{T} \mathbb{D}\vc{v}^L \qquad \text{ in } \Omega_L,
	\end{equation}
	\begin{equation}
		\label{2...}
		p^G - P^G = \frac{c^G}{T}\; \div \vc{v}^G, \qquad	\mathbb{S}^G = - p^G \mathbb{I}+ \frac{c^G}{T} \div \vc{v}^G \;\mathbb{I} + \frac{b^G}{T} \mathbb{D}\vc{v}^G \quad \text{ in } \Omega_G,
	\end{equation}
	and the Fourier's law,
	\begin{equation}
		\label{28}
		\vc{J}^L_q = -\frac{a^L}{T^2} \;\nabla_x T, \qquad \vc{J}^S_q = -\frac{a^S}{T^2} \;\nabla_x T, \qquad \vc{J}^G_q = -\frac{a^G}{T^2} \;\nabla_x T.
	\end{equation}
	Similarly, from the interface equations (\ref{25.})-(\ref{25...}), the following relations are deduced, neglecting the cross-diffusion effects:
	\begin{align}
		\vc{J}^s_{q,i} &= -\frac{a^s_i}{T^2} \;\nabla_\tau T \qquad i=1,2,3,	\label{28.} \\
		\mathbb{S}^s_i + P^s_i  \left(\mathbb{I} - \vc{n}_i \otimes \vc{n}_i\right) &= \frac{b^s_i}{T} \;\mathbb{D}\vc{v}^s_i \qquad i=1,2,3, \label{29} \\
		p^s_i - P^s_i &= \frac{c^s_i}{T}\left( \mathrm{div}_\Sigma \vc{v}^s_{i,\tau} +  \kappa_\Sigma v^s_{i,n}\right) \qquad i=1,2,3, \label{28..} \\
		\varrho^L(v^L - v^s_i)_n &=k_i T \Big[ \left( \mu^s_i - \mu^L\right) + \frac{1}{\varrho^L } \vc{n}_i\cdot \left( \mathbb{S}^L+p^L\mathbb{I}\right) \cdot\vc{n}_i\Big] \qquad i=1,2, \label{35} \\
		\frac{1}{2} \left[ \left( \mathbb{S}^L + \mathbb{S}^\chi\right) \cdot\vc{n}\right] _{\tau} &= \frac{\beta_i}{T} \left( \vc{v}^L - \vc{v}^\chi\right) _\tau \qquad \text{ on } \Sigma_i, \ (i,\xi) \in \{(1,G), (2,S)\}, \label{32.} 
	\end{align}
	\begin{align}
		\vc{v}^s_{i,\tau} = \frac{1}{2}\left( \vc{v}^L + \vc{v}^\chi\right) _\tau &-\alpha_i T \left[ \left( \mathbb{S}^L - \mathbb{S}^\chi\right) \cdot\vc{n}\right] _{\tau} \quad \text{ on } \Sigma_i, \ (i,\xi) \in \{(1,G), (2,S)\}, \label{32..} \\
		\frac{1}{2} \left[ \left( \mathbb{S}^G + \mathbb{S}^S\right) \cdot\vc{n}_3\right] _{\tau} &= \frac{\beta_3}{T} \left( \vc{v}^G - \vc{v}^S\right) _\tau \qquad \text{ on } \Sigma_3, \\
		\vc{v}^s_{3,\tau} &= \frac{1}{2}\left( \vc{v}^G + \vc{v}^S\right) _\tau -\alpha_3 T \left[ \left( \mathbb{S}^G - \mathbb{S}^S\right) \cdot\vc{n}_3\right] _{\tau} \qquad \text{on } \Sigma_3.
	\end{align}
	Here, the constants $a^\xi, b^\xi, c^\xi, a^s_i, b^s_i, c^s_i, k_i, \alpha_i, \beta_i $ are non-negative.
	
	\subsubsection{Further constitutive assumptions}
	
	In this section, we make some further assumptions that lead to some simplifications in the model.
	
	First of all, we assume that the ratio of the viscosities of the gas and the liquid is of the same order as the ratio of the corresponding densities\footnote{cf. https://www.engineersedge.com/physics/viscosity\_of\_air\_dynamic\_and\_kinematic\_14483.htm, https://www.accessengineeringlibrary.com/content/book/9780070471788/back-matter/appendix4}, that is,
	\begin{equation*}
		\frac{b^G}{b^L} \approx \Big(\frac{\varrho^G}{\varrho^L}\Big)^{1/2} \ll 1.
	\end{equation*}
	Hence, in (\ref{2...}), we can assume that $b^G=0$. Thus, from the momentum equation (\ref{pressure_gas}), one obtains, under the assumption that $c^G=0$ (which is the case for an ideal gas with $\varrho^G=0$),
	\begin{equation}
		\label{gas_pressure}
		p^G = \text{constant}.
	\end{equation}
	
	\paragraph{Interface stress tensors}
	The two-dimensional interface can be regarded as hydrodynamically ideal (cf. \cite[pp. 195]{shikh_bk}), i.e. there is no contribution of the viscous forces and thus the stress tensor is given by the pressure only
	. Therefore, equation (\ref{29}) reduces to,
	\begin{equation}
		\label{38}
		\mathbb{S}^s_i = -p^s_i (\mathbb{I} - \vc{n}_i\otimes \vc{n}_i), \qquad i=1,2,3,
	\end{equation}
	which implies
	\begin{equation}
		\label{2.}
		\mathrm{div}_\Sigma \mathbb{S}^s_i = -\nabla_\tau p^s_i - p^s_i \kappa_\Sigma \vc{n}_i.
	\end{equation}
	Here we assume that $b^s_i=0$ in (\ref{29}) and $c^s_i=0$ in (\ref{28..}), hence $P^s_i = p^s_i$. 
	The surface pressure $p^s_i$ is defined as the 'surface tension' taken with the negative sign.
	
	\paragraph{Interface mass balance equation}	In (\ref{35}), the last term on the RHS (hydrodynamic contribution) can be neglected compared to the first term on the RHS (chemical potential) (cf. \cite[pp. 196]{shikh_bk}). Also, in equilibrium, the following must hold,
	\begin{equation}
		\label{chem_equi}
		\mu^s_i(\varrho^s_{ie}, T) = \mu^L (\varrho^L,T), \quad i=1,2,
	\end{equation}
	where $ \varrho^s_{ie}$ is the equilibrium surface density.
	Further, assuming that $\mu^s_i$ can be expanded in a Taylor series around $\varrho^s_{ie}$, i.e.
	\begin{equation}
		\label{surface_density}
		\mu^s_i = \mu^s_i(\varrho^s_{ie}) + \frac{\partial \mu^s_i}{\partial \varrho^s_i}\biggr\rvert_{\varrho^s_{ie}}\left( \varrho^s_i - \varrho^s_{ie}\right) ,
	\end{equation}
	(\ref{35}) can be written as,
	\begin{equation}
		\label{36}
		\varrho^L(v^L_n - v^s_{i,n}) = \frac{\varrho^s_i - \varrho^s_{ie}}{\tau_i}, \quad i=1,2,
	\end{equation}
	where $\tau_i$
	is a relaxation time. Note that although the same notation is used to denote the relaxation time and the tangent vectors, we use $\tau$ as indices over a vector to indicate tangential component, so there should not be any ambiguity.
	Using (\ref{36}), we can now simplify the surface mass balance equation (\ref{37}) as
	\begin{equation}
		\label{42}
		\partial_t^\Sigma \varrho^s_i + \mathrm{div}_\Sigma \left( \varrho^s_i \vc{v}^s_{i,\tau}\right)  + \kappa_\Sigma \varrho^s_i v^s_{i,n} =-\frac{\varrho^s_i - \varrho^s_{ie}}{\tau_i}, \qquad i=1,2.
	\end{equation}
	
	\paragraph{Interface momentum equation}
	For the momentum balance equations (\ref{17...})-(\ref{17..}), the surface convective momentum $\partial_t^\Sigma \big( \varrho^s_i \vc{v}^s_i\big)  + \mathrm{div}_\Sigma \left( \varrho^s_i \vc{v}^s_i \otimes \vc{v}^s_{i,\tau}\right) + \kappa_\Sigma \varrho^s_i \vc{v}^s_i v^s_{i,n} $ can be neglected compared to the capillary effect $\mathrm{div}_\Sigma \mathbb{S}^s_i$ (cf. \cite[pp. 197]{shikh_bk}). 
	Also the contribution of the term $\varrho^L\vc{v}^L(v^L_n - v^s_{i,n})$ on the left hand side can be neglected (cf. \cite[pp. 198]{shikh_bk})
	, leading to the form,
	\begin{equation}
		\label{39}
		\mathrm{div}_\Sigma \mathbb{S}^s_i + \big( \mathbb{S}^L - \mathbb{S}^\xi\big) \cdot\vc{n}_i = 0 \qquad \text{ on } \Sigma_i, \ (i,\xi) \in \{(1,G), (2,S)\},
	\end{equation}
	\begin{equation*}
		\mathrm{div}_\Sigma \mathbb{S}^s_3 + \left( \mathbb{S}^G - \mathbb{S}^S\right) \cdot\vc{n}_3 = 0 \quad \text{ on } \Sigma_3,
	\end{equation*}
	which can be splitted into the normal and tangential component, invoking the expression for the liquid stress tensor (\ref{2..}), relation (\ref{2.}) and the incompressibility assumption of the liquid, as
	\begin{equation}
		\label{43}
		\left( p^\xi - p^L\right) + \frac{b^L}{T} \;\vc{n}_i\cdot \left( \mathbb{D}\vc{v}^L\right)\cdot\vc{n}_i + p^s_i \kappa_\Sigma = 0 \qquad  (i,\xi) \in \{(1,G), (2,S)\},
	\end{equation}
	and
	\begin{equation}
		\label{40}
		\frac{b^L}{T} \left[ \left( \mathbb{D}\vc{v}^L\right)\cdot\vc{n}_i\right]_\tau - \nabla_{\tau} p^s_i= 0 \qquad i=1,2.
	\end{equation}
	Equation (\ref{40}) is the well-known Marangoni effect which is the tangential forces caused due to the surface tension gradient.
	
	\paragraph{Velocity equation}	Next, eliminating $\vc{v}^\chi$ from (\ref{32.}) and (\ref{32..}), we get with the help of equation (\ref{39}), substituting $\mathbb{S}^L$,
	\begin{equation}
		\label{44}
		T \left( 1+ 4\alpha_1 \beta_1 \right) \nabla_{\tau} p^s_1 + 4\beta_1 \left( \vc{v}^s_{1,\tau} - \vc{v}^L_\tau\right) =0 \qquad \text{ on } \Sigma_1.
	\end{equation}
	
	\paragraph{Barotropic condition}	Lastly, assuming that the interface formation process can be modelled as a barotropic process, a linear equation of state is considered which reflects the changes of the surface pressure due to the compression or extreme rarefaction of the surface phase,
	\begin{equation}
		\label{45}
		p^s_i = \gamma_i \left( \varrho^s _{i,0} - \varrho^s_i\right), \qquad i=1,2 .
	\end{equation}
	The parameter $\gamma_i$ is a constant describing compressibility of the fluid at the interface and $\varrho^s_{i,0}$ is the surface density corresponding to zero surface pressure.
	
	\paragraph{Liquid-solid interface} 	In a similar manner, we can write the boundary conditions on the liquid-solid interface $\Sigma_2$ as well, such as (\ref{42}), (\ref{36}), (\ref{45}), together with (\ref{normal_velocity}). The different condition is a generalized Navier slip condition (instead of (\ref{44})) which one obtains from (\ref{32.}) and (\ref{39}) eliminating $\mathbb{S}^S$,
	\begin{equation}
		\label{23}
		\frac{b^L}{T}\left[ \left( \mathbb{D}\vc{v}^L\right)\cdot\vc{n}_2\right]_\tau - \frac{1}{2}\nabla_{\tau} p^s_2= \frac{\beta_2}{T} \left( \vc{v}^L_\tau - \vc{v}^S_\tau\right).
	\end{equation}
	together with (\ref{32..}).	This condition shows that there can be apparent slip on the solid substrate.
	
	\paragraph{Gas-solid interface} 
	
	Also, we assume that the density and the chemical potential of the gas-solid interface do not make any contribution to the liquid and contact line dynamics, therefore,
	\begin{equation}
		\label{gas_solid_surface}
		\varrho^s_3=0, \quad \mu^s_3 =0.
	\end{equation}
	However, the contribution of the surface tension $p^s_3$ remains.
	
	\paragraph{Equations at the contact line}	Finally, for the equations at the contact line,
	taking into account the above assumptions, the mass equation (\ref{mass_cont}) reduces to,
	\begin{equation}
		\label{mass_cont1}
		\varrho^s_1 \left( \vc{v}^s_{1,\tau} - \vc{v}^c\right)\cdot \vc{\nu}_1 +\varrho^s_2 \left( \vc{v}^s_{2,\tau} - \vc{v}^c\right)\cdot \vc{\nu}_2 =0 \quad \text{ on } C.
	\end{equation}
	
	Next we re-write the momentum balance equation (\ref{64}). 
	Let $\theta_i$ be the angles at which the interfaces $\Sigma_i$ intersect the contact line. Therefore,
	\begin{equation*}
		\cos \theta_1 = \vc{\nu}_1\cdot \vc{\nu}_2, \qquad \cos \theta_2 = \vc{\nu}_2\cdot \vc{\nu}_3, \qquad \cos \theta_3 = \vc{\nu}_1\cdot \vc{\nu}_3.
	\end{equation*}
	Then, as the convective momentum fluxes along the interfaces are negligible compared to the surface tensions, we can write (\ref{64}) as,
	\begin{equation*}
		0 = \sum_{i=1}^3 \mathbb{S}^s_i \cdot\nu_i = - \sum_{i=1}^3 p^s_i (\mathbb{I} - \vc{n}_i\otimes \vc{n}_i) \cdot \vc{\nu}_i,
	\end{equation*}
	which gives taking projecting on $\vc{\nu}_2$ direction,
	\begin{equation*}
		0=	\sum_{i=1}^3 \vc{\nu}_2 \cdot \mathbb{S}^s_i \cdot \vc{\nu}_i = -p^s_1 \cos \theta_1 -p^s_2 -p^s_3 \cos \theta_2.
	\end{equation*}
	In the particular situation when $\theta_2 = \pi$, the above relation yields the Young equation involving the contact angle $\theta$ between the liquid-gas and liquid-solid interface,
	\begin{equation}
		\label{young}
		p^s_1 \cos \theta=p^s_3-p^s_2.
	\end{equation}
	
	For the entropy equation at the contact line (\ref{CL_entropy}), using the energy equation (\ref{energy_cont}) together with the identity for heat flux (\ref{fluxes}) 
	and the Gibbs relation (\ref{20.}), it becomes,
	\begin{equation*}
		\begin{aligned}
			\sum_{i=1}^3\frac{1}{T}\Big[ (\vc{v}^s_{i,\tau} - \vc{v}^c)\cdot (\mathbb{S}^s_i \cdot \nu_i)\Big] -\frac{1}{T} \Big[ \varrho^s_i (\vc{v}^s_{i,\tau} - \vc{v}^c)\cdot\nu_i \Big( \mu^s_i + \frac{|\vc{v}^s_i - \vc{v}^c|^2}{2} - \frac{1}{\varrho^s_i} p^s_i\Big) \Big] =0.
		\end{aligned}
	\end{equation*}
	Further incorporating the interface stress relation (\ref{38}), the above equation reduces to,
	\begin{equation}
		\label{entropy_cont}
		\sum_{i=1}^3 \frac{1}{T} \Big[ \varrho^s_i (\vc{v}^s_{i,\tau} - \vc{v}^c)\cdot\nu_i \Big( \mu^s_i + \frac{|\vc{v}^s_i - \vc{v}^c|^2}{2} \Big) \Big] = 0.
	\end{equation}
	The square term can be neglected compared to the chemical potential in the above relation, due to the similar argument as used before to deduce (\ref{36}) from (\ref{35}) (cf. \cite[pp. 196]{shikh_bk}).
	Moreover, employing the assumptions concerning the gas-solid interface (\ref{gas_solid_surface}) and the mass balance condition at the triple junction (\ref{mass_cont1}), relation (\ref{entropy_cont}) yields,
	\begin{equation}
		\label{chem_jump}
		\varrho^s_1 \;\mu^s_1 (\vc{v}^s_{1,\tau} - \vc{v}^c)\cdot\nu_1 + \varrho^s_2 \;\mu^s_2 (\vc{v}^s_{2,\tau} - \vc{v}^c)\cdot\nu_2 =0.
	\end{equation}
	But we also have the following expression for the chemical potential in terms of the density (cf. (\ref{chem_equi}), (\ref{surface_density})),
	\begin{equation*}
		\mu^s_i = \mu^s_{i,e} + m_i \left( \varrho^s_i - \varrho^s_{ie}\right) = \mu^L_e+ m_i \left( \varrho^s_i - \varrho^s_{ie}\right).
	\end{equation*}
	Therefore, we obtain from the jump condition for the chemical potential (\ref{chem_jump}), with the assumption that $m_2\ne 0$ (i.e. the chemical potential varies from its equilibrium value),
	\begin{equation*}
		\frac{m_1}{m_2} \varrho^s_1 \left( \varrho^s_1 - \varrho^s_{1,e}\right) (\vc{v}^s_{1,\tau} - \vc{v}^c)\cdot\nu_1 + \varrho^s_2 \left( \varrho^s_2 - \varrho^s_{2e}\right)(\vc{v}^s_{2,\tau} - \vc{v}^c)\cdot\nu_2 =0.
	\end{equation*}
	But we can further simplify the above relation with the help of the mass equation (\ref{mass_cont1}) to,
	\begin{equation}
		\label{missing_cond}
		m\left( \varrho^s_1 - \varrho^s_{1,e}\right) - \left( \varrho^s_2 - \varrho^s_{2e}\right)=0,
	\end{equation}
	where $m=	\frac{m_1}{m_2} $ and	under the assumption that $\varrho^s_1\ne 0$ and $(\vc{v}^s_{1,\tau} - \vc{v}^c)\cdot\nu_1 \ne 0$.
	
	\subsubsection{Isothermal limit}
	We investigate now the situation in the case of an isothermal process more carefully.
	In particular, the goal is to show that in this regime, the temperature is almost constant in the liquid and consequently, the liquid energy equation and liquid-solid and the liquid-gas surface energy equations do not contribute to the model any more.
	
	First, let us consider that the temperature of the liquid can be expressed as
	\begin{equation*}
		T= T_0\Big( 1 + \frac{1}{	\tilde{a}^L}\; \theta^L\Big) \qquad \text{ in } \Omega_L,
	\end{equation*}
	where $T_0$ is a constant temperature, $\theta^L$ is a small deviation from the constant temperature, and $\tilde{a}^L$ is a non-dimensional quantity given by,
	\begin{equation}
		\label{a_L}
		\begin{aligned}
			\tilde{a}^L = \min \left\lbrace 
			\frac{[a^L][dT] }{[b^L] [T][\vc{v}^L]^2},
			\; \frac{[a^L][dT]}{L [T]^3 [\varrho^L] [\vc{v}^L] [c^L_m]}, \; \frac{[a^L] L}{[a^s_2]}, \; \frac{[a^L] L}{[a^s_1]}, \; \frac{[a^L][dT]}{L[T]^2[\varrho^L][\vc{v}^L]^3}, \; \frac{[a^L]}{[a^G]}
			\right\rbrace.
		\end{aligned}
	\end{equation}
	We recall that $a^L$ is the heat conductivity of the liquid, appearing in (\ref{28}), $L$ is the characteristic length scale of the underlying system, $c^L_m$ is the specific heat capacity of the liquid (cf. (\ref{internal_energy_temperature})), and $[\cdot]$ denotes the dimension of a physical quantity.
	The first quantity of the RHS in (\ref{a_L}) represents the ratio of the effect of the thermal conductivity to the source term due to the viscous forces in the equation (\ref{22}), while the second term is the ratio of the thermal conductivity to the transport of heat due to convection. The other terms in the above expression are chosen in such a way in order to make the contribution of the heat flux from the liquid in the surface energy equation (\ref{11.}) much larger, compared to all the other terms (as explained below).
	
	In obtaining (\ref{a_L}), the following facts have been used. First of all, $\varrho^s_i\simeq \varrho^L \ell$ where $\ell$ denotes the thickness of the interface and $\ell \ll L$; The surface velocity is of the same order as the bulk velocity and also $[\vc{v}^s_i]=[\vc{v}^L]$; Similarly, the heat capacities of the interface and of the bulk is of the same order and $[c^L_m] = [c^s_{m,i}]$; And importantly, the temperature difference between two bulk phases is small compared to the temperature in the liquid $dT \ll T$. This last assumption is necessary in order to obtain isothermal limit. Otherwise there might be large heat flux which could make the solid temperature far from constant.
	
	With the above definition, if we assume that $\tilde{a}^L\gg 1$, it means that conduction is much bigger than both the convection and the source due to the viscous force in the internal energy equation of the liquid (\ref{22}).			Therefore, all the other terms in (\ref{22}) become negligible compared to the heat flux and in the limit $\tilde{a}^L \to \infty$, we obtain a heat equation,
	\begin{equation}
		\label{heat}
		-\div \Big( \frac{1}{T^2} \nabla_{\mathbf{x}} T\Big) = 0 \qquad \text{ in } \Omega_L.
	\end{equation}
	
	On the liquid-solid interface $\Sigma_2$,		 
	combining the mass density equation (\ref{37}) with the internal energy equation (\ref{11.}) for $(i,\xi)=(2,S)$, we obtain, taking into account the specific forms of the heat fluxes (\ref{28}), (\ref{28.}), and also the assumption (\ref{temperature_continuity}),
	\begin{align}
		& \quad \varrho^s_2\left( \partial^\Sigma_t u^s_2 + \vc{v}^s_{2,\tau} \cdot \nabla_{\tau} u^s_2\right) + \varrho^L \left( u^L - u^s_2\right) \left(v^L_n - v^s_{2,n}\right) - a^s_2\; \mathrm{div}_\Sigma \Big( \frac{1}{T^2} \nabla_{\tau}T\Big) \nonumber\\
		& = \vm{S}^s_2 : \nabla_{\tau}\vc{v}^s_2 +\vc{n}_2\cdot \vm{S}^L \cdot\left( \vc{v}^L - \vc{v}^s_2\right) -\vc{n}_2 \cdot \vm{S}^S \cdot\left( \vc{v}^S - \vc{v}^s_2\right) +\frac{a^L}{T^2} \frac{\partial T}{\partial \vc{n}}\Bigr|_{\text{liquid}}-\frac{ a^S}{T^2} \frac{\partial T}{\partial \vc{n}}\Bigr|_{\text{solid}}
		\nonumber \\
		& \hspace{0.3cm}- \frac{1}{2}\varrho^L (v^L_n - v^s_{2,n}) (\vc{v}^s_2 -\vc{v}^L)^2 \hspace{5cm} \text{ on } \Sigma_2. \label{19}
	\end{align}
	Note that, although the temperature is continuous across the interfaces, there might be a jump for the normal derivative, i.e. $\frac{\partial T}{\partial \vc{n}}\Bigr|_{\text{liquid}} \ne \frac{\partial T}{\partial \vc{n}}\Bigr|_{\text{solid}}$ in general. Here the subscript means the quantity is evaluated at the interface $\Sigma_2$ in a limiting sense, with values approaching from the liquid or the solid, i.e. $\frac{\partial T}{\partial \vc{n}}\Bigr|_{\text{liquid}} = \underset{x\to x_0, x\in \Omega_L}{\lim} \frac{\partial T}{\partial \vc{n}}(x), x_0\in \Sigma_2$. In the limit $\tilde{a}^L\to \infty$, due to the definition (\ref{a_L}), 
	the above energy equation reduces to nothing but a continuity condition,
	\begin{equation}
		\label{LS_BC}
		a^L \;\frac{\partial T}{\partial \vc{n}}\Bigr|_{\text{liquid}}= a^S \;\frac{\partial T}{\partial \vc{n}}\Bigr|_{\text{solid}} \qquad \text{ on $\Sigma_2$ }.
	\end{equation}
	
	For the liquid-gas interface $\Sigma_1$, we can perform similar analysis as above. At this point, note that although we obtained that the pressure of the gas is constant, it does not necessarily imply that the internal energy (and hence, the temperature) is constant too, due to the ideal gas law. One may only have the following relation,
	$$
	\varrho^G u^G = \frac{3}{2} R^G \varrho^G T = \frac{3}{2} p^G \approx p^L \ll \varrho^L u^L,
	$$
	where $R^G$ is the ideal gas constant.
	Plugging in the relations for heat fluxes (\ref{28}),  (\ref{28.}), the energy equation (\ref{11.}) for $(i,\xi)=(1,G)$ transforms to, together with the mass equation (\ref{37}),
	\begin{equation}
		\label{liquid_gas}
		\begin{aligned}
			& \varrho^s_1\left( \partial^\Sigma_t u^s_1 + \vc{v}^s_{1,\tau} \cdot \nabla_{\tau} u^s_1\right) + \varrho^L \left( u^L - u^s_1\right) \left(v^L_n - v^s_{1,n}\right)- a^s_1\; \mathrm{div}_\Sigma \Big( \frac{1}{T^2} \nabla_{\tau}T\Big) \\
			= & \; \vm{S}^s_1 : \nabla_{\tau}\vc{v}^s_1 +\vc{n}_1\cdot \mathbb{S}^L\cdot(\vc{v}^L - \vc{v}^s_1) - p^G(v^G_n - v^s_{i,n})+\frac{a^L}{T^2} \frac{\partial T}{\partial \vc{n}}\Bigr|_{\text{liquid}}-\frac{ a^G}{T^2} \frac{\partial T}{\partial \vc{n}}\Bigr|_{\text{gas}}\\
			& \quad - \frac{1}{2}\varrho^L (v^L_n - v^s_{1,n}) (\vc{v}^L-\vc{v}^s_1)^2 \hspace{6cm} \text{ on } \Sigma_1.
		\end{aligned}
	\end{equation}
	Thus, in the limit $\tilde{a}^L \to \infty$, because of the definition (\ref{a_L}), all the other terms in the equation (\ref{liquid_gas}) can be neglected compared to the temperature gradient from the liquid and
	we obtain a homogeneous Neumann boundary condition,
	\begin{equation}
		\label{LG_BC}
		\frac{\partial T}{\partial \vc{n}}\Bigr|_{\text{liquid}}=0 \qquad \text{on $\Sigma_1$}.
	\end{equation} 
	
	Next, note that the solid energy equation (\ref{22}) becomes, since $\div \vc{v}^S=0, \mathbb{D}\vc{v}^S =0$ due to the structure (\ref{solid_velocity}), and taking into account the relation between the changes of internal energy and the temperature by $\delta u^S = c^S_m \; \delta T$ with $c^S_m$ being the specific heat capacity (i.e. per unit mass),
	\begin{equation}
		\label{solid_energy}
		c^S_m \varrho^S \left( \partial_t T + \vc{v}^S\cdot \nabla_x T \right)  =  a^S\div\Big( \frac{1}{T^2} \;\nabla_x T \Big) \quad \mbox{ in } \Omega_S.
	\end{equation}
	Defining the non-dimensional parameter as before,
	\begin{equation}
		\label{a_S}
		\tilde{a}^S =
		\min \left\lbrace \frac{[a^S]}{L [T]^2 [\varrho^S] [\vc{v}^S] [c^S_m]}, \;
		\frac{[a^S][\delta T] }{L [T]^2[p^G][\vc{v}^G]},
		\; \frac{[a^S] L}{[a^s_3] }, \; \frac{[a^S]}{[a^G]}
		\right\rbrace,
	\end{equation}
	and letting $\tilde{a}^S \to \infty$, i.e. assuming that the effect of the conduction is much bigger than the convection, the solid energy equation (\ref{solid_energy}) further reduces to a heat equation,
	\begin{equation}
		\label{solid_heat}
		\div\Big( \frac{1}{T^2} \;\nabla_x T \Big)=0 \quad \mbox{ in } \Omega_S.
	\end{equation}
	
	Similarly, in the limit $\tilde{a}^S \to \infty$, the gas-solid energy equation (\ref{GS_energy}) becomes a Neumann condition, since we essentially make the contributions of all the other terms negligible compared to the heat flux from the solid with the above definition (\ref{a_S}) of $\tilde{a}^S$, 
	\begin{equation}
		\label{GS_BC}
		\frac{\partial T}{\partial \vc{n}}\Bigr|_{\text{solid}}=0 \qquad \text{on $\Sigma_3$}.
	\end{equation}
	
	From the above discussion, we obtain that the heat equations (\ref{heat}), (\ref{solid_heat}) in the bulk together with the Neumann conditions (\ref{LG_BC}), (\ref{LS_BC}) and (\ref{GS_BC}) imply that the temperature in the liquid and the solid is constant $T_0$. Thus, the energy equations in the bulk and at the interfaces are completely decoupled from the main contact line problem and do not appear in the final model.
	
	\paragraph{Missing condition:}
	
	The energy equation at the contact line (\ref{missing_cond}), however, remains as it is, even in the isothermal limit. This third boundary condition has not been taken into account in the original model in \cite{shikh93} due to the assumption of isothermal flows (although mentioned in the book \cite[p 210]{shikh_bk} very briefly, claiming that the quantities such as "surface energy" are artificial concepts in a macroscopic model, and introducing one would result in an error; See also pp. 219). As mentioned in the introduction, the missing condition was first pointed out in \cite{bil06} and was shown to be necessary for well-posedness of the contact line problem (cf. \cite[discussion under "final remarks"]{bothe}). 
	However, our derived missing condition is somewhat different from the one stated in \cite[equation (2.23) or (2.24)]{bil06}, but same as in \cite[equation (66)]{bothe}.
	
	\subsection{Summary of the model}
	To summarize, we recall below the complete system of the moving contact line, for an incompressible isothermal liquid. Henceforth, we re-name some of the constants, such as $\frac{b^L}{T} \to 2\mu$ for the dynamic viscosity of the fluid (which must not be confused with the chemical potential, which in any case is not used explicitly furthermore) and $\frac{\beta_i}{T} \to \beta_i, \alpha_i T \to \alpha_i$ to be consistent with the standard notations. \\

	\noindent\textbf{equation of the liquid:} Equations (\ref{incom}) and (\ref{NS}), together with (\ref{2..}) give,	
	
	\begin{equation}
		\label{1}
		\div \;\vc{v}^L =0; \quad \quad \varrho^L\left( \partial_t \vc{v}^L + \vc{v}^L\cdot \nabla_{\mathbf{x}}\vc{v}^L\right) = \mu\Delta \vc{v}^L - \nabla_{\mathbf{x}} p^L;
	\end{equation}
	
	\noindent \textbf{equation at the liquid-gas interface:} It follows from equations (\ref{46}), (\ref{42}), (\ref{36}), (\ref{43}), (\ref{40}), (\ref{44}) and (\ref{45}) that
	\begin{subequations}
		\begin{align}
			V_n &= \vc{v}^s_1\cdot \vc{n}_1, \label{2}\\
			\partial_t^\Sigma \varrho^s_1 + \mathrm{div}_\Sigma \left( \varrho^s_1 \vc{v}^s_{1,\tau}\right)  + \kappa_\Sigma \varrho^s_1 v^s_{1,n} &=-\frac{\varrho^s_1 - \varrho^s_{1e}}{\tau_1}, \label{3}
		\end{align}
		\begin{align}
			\varrho^L(v^L_n - v^s_{1,n}) &= \frac{\varrho^s_1 - \varrho^s_{1e}}{\tau_1}, \label{4}\\
			\left( p^G - p^L\right) + 2\mu\left[ \mathbb{D}\vc{v}^L \cdot\vc{n}_1\right] \cdot \vc{n}_1 + p^s_1 \kappa_\Sigma &= 0, \label{5}\\
			2 \mu\left[ \mathbb{D}\vc{v}^L\cdot\vc{n}_1\right]_\tau - \nabla_{\tau} p^s_1 &= 0, \label{6}\\
			-\left( 1+ 4\alpha_1 \beta_1 \right) \nabla_{\tau} p^s_1 &= 4\beta_1 \left( \vc{v}^s_{1,\tau} - \vc{v}^L_\tau\right), \label{7}\\
			p^s_1 &= \gamma_1 \left( \varrho^s_1 - \varrho^s_{1,0}\right), \label{8}
		\end{align}
	\end{subequations}
	
	\noindent \textbf{equation at the liquid-solid interface:} Similarly, equations (\ref{42}), (\ref{36}), (\ref{23}), (\ref{normal_velocity}), (\ref{32..}) and (\ref{45}) yield,
	\begin{subequations}
		\begin{align}
			\partial_t^\Sigma \varrho^s_2 + \mathrm{div}_\Sigma \left( \varrho^s_2 \vc{v}^s_{2,\tau}\right) + \kappa_\Sigma \varrho^s_2 v^s_{2,n} &=-\frac{\varrho^s_2 - \varrho^s_{2e}}{\tau_2}, \label{9}\\
			\varrho^L(v^L_n - v^s_{2,n}) &= \frac{\varrho^s_2 - \varrho^s_{2e}}{\tau_2}, \label{12}\\
			2 \mu\left[  \mathbb{D}\vc{v}^L\cdot\vc{n}_2\right]_\tau - \frac{1}{2}\nabla_{\tau} p^s_2 &= \beta_2 \left( \vc{v}^L_\tau - \vc{v}^S_\tau\right) , \label{13}\\
			\left( \vc{v}^s_2 - \vc{v}^S\right) \cdot \vc{n}_2 &=0, \label{14}\\
			\vc{v}^s_{2,\tau} &= \frac{1}{2}\left( \vc{v}^L + \vc{v}^S\right) _\tau -\alpha_2 \nabla_{\tau} p^s_2, \label{15}\\
			p^s_2 &= \gamma_2 \left( \varrho^s_2 - \varrho^s_{2,0}\right), \label{16}
		\end{align}
	\end{subequations}
	
	\noindent \textbf{equation at the triple junction:} Equations (\ref{mass_cont1}), (\ref{young}) and (\ref{missing_cond}) give,
	\begin{subequations}
		\begin{align}
			\varrho^s_1\left( \vc{v}^s_{1,\tau} - \vc{v}^c\right) \cdot \vc{\nu}_1 + \varrho^s_2 \left( \vc{v}^s_{2,\tau} - \vc{v}^c\right) \cdot \vc{\nu}_2 &=0, \label{17}\\
			p^s_1 \cos \theta &= p^s_3 - p^s_2, \label{18}\\
			m \left( \varrho^s_1 - \varrho^s_{1e}\right) - \left( \varrho^s_2 - \varrho^s_{2e}\right) &=0, \label{18.1}
		\end{align}
	\end{subequations}
	
	\section{Lubrication approximation}
	\label{S3}
	
	The general three-dimensional problem can be considerably simplified in the thin film approximation where the ratio of the characteristic length scales in the vertical and horizontal direction to the film is considered as a small parameter. Below we derive different asymptotic limits depending on different orders of magnitude for the slip length parameter and also for different models.
	
	\subsection{Thin film model for small slip length}
	\label{S2.1}
	
	First we consider the simple case of spreading of a viscous, incompressible droplet of liquid over a smooth solid surface where the displaced medium is an inviscid gas, in two dimension.

	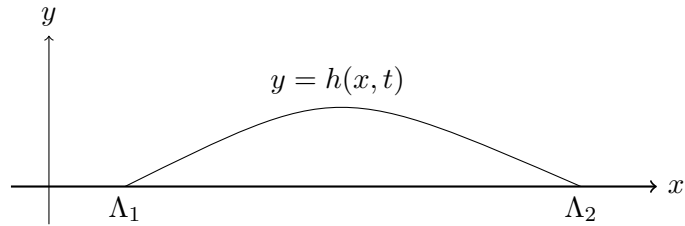
\begin{figure}[h]
		\centering
		\begin{tikzpicture}
			\draw[->,thick] (-0.5,0)--(8,0) node[right]{$x$};
			\draw[->] (0,-0.5)--(0,2) node[above]{$y$};
			\draw (1,0) .. controls (3.8,1.4) .. (7,0);
			\node at (3.8, 1.4) {$y=h(x,t)$};
			\node at (1,-0.3) {$\Lambda_1$};
			\node at (7, -0.3) {$\Lambda_2$};
		\end{tikzpicture}
		\caption{thin film model for droplet}
	\end{figure}
	
	The liquid-gas interface is given by $y = h(x,t)$ which is the free surface and the liquid-solid interface is $y=0$. The contact points where these interfaces meet are $s_i(t) =(\Lambda_i(t),0) $. We denote by $L$ the characteristic length scale of the droplet, $H$ the characteristic thickness of the film with $H\ll L$ and $ \varepsilon = \frac{H}{L}$. We are interested in obtaining the limiting model as $\varepsilon \to 0$.
	
	In Shikhmurzaev's model (\ref{1})-(\ref{18}), we use the notation $ \vc{v}^L= (u,v), \vc{v}^s_i= (u^s_{i},v^s_{i}), i=1,2$ and we assume $\vc{v}^S=0$, i.e. the fluid droplet is on a plane at rest. We omit the index for the pressure of the liquid and just denote it by $p$.
	Let us then consider the following scales,
	\begin{align*}
		&\overline{x} = \frac{x}{L}, && \overline{y} = \frac{y}{H}, && \overline{\Lambda}_i = \frac{1}{L} \Lambda_i, && \overline{h} = \frac{h}{H},\\[.3cm]
		&\overline{t} = \frac{\varepsilon^3 \sigma_{1e}}{L \mu} t, && \overline{u} = \frac{\mu}{\varepsilon^3\sigma_{1e}} u, && \overline{v} = \frac{\mu}{\varepsilon^4 \sigma_{1e}} v, &&
		\overline{p} =\frac{ L}{\varepsilon \sigma_{1e}} p,\\[.3cm]
		&  \overline{u}^s_{i} = \frac{\mu}{\varepsilon^3\sigma_{1e}} u^s_{i}, && \overline{v}^s_{i} = \frac{\mu}{\varepsilon^4 \sigma_{1e}}v^s_{i}, && (\varrho^s_i - \varrho^s_{ie}) = \frac{\varepsilon^2 \sigma_{1e}}{\gamma_1}\zeta_i, \quad
		&& \theta_d =\varepsilon \partial_{\overline{x}} \overline{h}.
	\end{align*}
	Here
	$ \sigma_{1e} := - p^s_1 (\varrho^s_{1e}) = - \gamma_1 (\varrho^s_{1e} - \varrho^s_{1,0}) >0$.
	The time scale is chosen in such a way in order to make the characteristic time of spreading of the droplet of order one. With the time and length scales fixed, it is then standard to re-scale the velocities accordingly. The pressure is scaled as in the standard thin film approximation, in order not to make the first component of the Navier-Stokes equation trivial 
	(cf. $(\ref{20})_3$). Also, the scaling for difference of the surface densities is made due to the relation with the surface tensions (\ref{8}) or, (\ref{16}).
	
	\subsubsection{Reduced model}
	
	We show below that equations (\ref{1})-(\ref{18.1}) reduce, to leading order terms in $\varepsilon$, to the following system:\\
	in $(\overline{x}, \overline{y}) \in (\overline{\Lambda}_1, \overline{\Lambda}_2)\times (0,\overline{h})$:
	\begin{equation}
		\label{20}
		\partial_{\overline{x}} \overline{u} + \partial_{\overline{y}} \overline{v} =0, \quad 
		\partial^2_{\overline{y}} \overline{u} = \partial_{\overline{x}} \overline{p}
		, \quad \partial_{\overline{y}} \overline{p}=0;
	\end{equation}
	on $\overline{y} = \overline{h}$:
	\begin{subequations}
		\begin{align}
			\partial_{\overline{t}} \overline{h} + \overline{u}^s_{1} \;\partial_{\overline{x}} \overline{h} &= \overline{v}^s_{1}, \label{21.1}\\
			\partial_{\overline{x}} \overline{u}^s_{1} &= -\frac{1}{\lambda}_1 \zeta_1, \label{21.2}\\
			\overline{u}^s_{1} \;\partial_{\overline{x}} \overline{h} - \overline{v}^s_{1} &= \overline{u} \;\partial_{\overline{x}} \overline{h} - \overline{v}, \label{21.3}\\
			\overline{p} &= - \;\partial^2_{\overline{x}} \overline{h}, \label{22.1}\\
			\partial_{\overline{y}} \overline{u} &= -\partial_{\overline{x}} \zeta_1, \label{22.2}\\
			\left( \overline{u}^s_{1} - \overline{u} \right) &=- d_1 \; \partial_{\overline{x}} \zeta_1; \label{22.3}
		\end{align}
	\end{subequations}
	on $\overline{y}=0$:
	\begin{subequations}
		\begin{align}
			\partial_{\overline{x}} \overline{u}^s_{2} &=-\frac{1}{\lambda}_2 \zeta_2, \label{23.1}\\
			\overline{v} &= \overline{v}^s_{2}, \label{23.2}\\
			\partial_{\overline{y}} \overline{u} &= \frac{1}{2}a\; \partial_{\overline{x}} \zeta_2+\beta\; \overline{u} , \label{31.1}\\
			\overline{v}^s_{2} &=0, \label{31.3}\\
			\overline{u}^s_{2} &= \frac{1}{2} \overline{u} - d_2 \;\partial_{\overline{x}} \zeta_2; \label{31.2}
		\end{align}
	\end{subequations}
	at $(\overline{x},\overline{y})=(\overline{\Lambda}_i,0), i=1,2$:
	\begin{subequations}
		\begin{align}
			\varrho^s_{1e} \left(\overline{u}^s_{1} - \partial_{\overline{t}}\overline{\Lambda}_i\right) + \varrho^s_{2e} \left(\overline{u}^s_{2} -\partial_{\overline{t}}\overline{\Lambda}_i\right) &=0 ,\label{44.1}\\
			\zeta_1 +a\;\zeta_2+ \frac{1}{2} |\partial_{\overline{x}} \overline{h}|^2 &= b, \label{44.3}\\
			m \;\zeta_1 &= \zeta_2.
			\label{44.4}
		\end{align}
	\end{subequations}
	In the above equations, $d_1 = \frac{ \mu}{L}\frac{(1+ 4\alpha_1 \beta_1)}{4\varepsilon\beta_1},
	\frac{1}{\lambda_1} = \frac{L \mu}{\varepsilon \varrho^s_{1e}\gamma_1 \tau_1},
	\frac{1}{\lambda_2} = \frac{L \mu}{\varepsilon \varrho^s_{2e}\gamma_1 \tau_2},
	a= \frac{\gamma_2}{\gamma_1},
	\beta = \varepsilon \beta_2 \frac{L}{\mu},
	d_2 = \alpha_2\frac{\gamma_2}{\varepsilon \gamma_1} \frac{\mu}{L},
	b =- \frac{\sigma_{3e}}{\sigma_{1e}} \overline{ \sigma}_3 $
	and all these non-negative constants 
	are $O(1)$ and dimensionless. The derivation of the above system is sketched below.
	
	To derive (\ref{20})-(\ref{44.4}), let us first compute some related terms in rescaled variables. The free surface, i.e. the fluid-gas interface is now given by $\overline{y} = \overline{h}(\overline{x},\overline{t})$, and the unit normal vector, outward with respect to the gas domain and the corresponding tangent vector are given by
	\begin{equation}
		\label{normal_tangent_vector}
		\vc{n}_1 = \frac{\left( \varepsilon \partial_{\overline{x}} \overline{h}, -1\right) }{(1+ \varepsilon^2 |\partial_{\overline{x}} \overline{h}|^2)^{1/2}}, \qquad \vc{\tau} = \frac{\left(1,  \varepsilon \partial_{\overline{x}} \overline{h}\right) }{(1+ \varepsilon^2 |\partial_{\overline{x}} \overline{h}|^2)^{1/2}},
	\end{equation}
	and the curvature by
	\begin{equation}
		\label{curvature}
		\kappa_\Sigma = \frac{\varepsilon \;\partial^2_{\overline{x}} \overline{h}}{L\left( 1+ \varepsilon^2 |\partial_{\overline{x}} \overline{h}|^2\right) ^{3/2}}.
	\end{equation}
	We compute the rate of strain tensor
	\begin{equation*}
		\renewcommand{\arraystretch}{1.5}
		\mathbb{D}\vc{v}^L = 
		\frac{1}{2} \begin{bmatrix}
			2\partial_{x} u & \partial_{y} u + \partial_{x} v\\
			\partial_{y} u + \partial_{x} v & 2 \partial_{y} v
		\end{bmatrix}
		=\frac{1}{2} \frac{\varepsilon^3 \sigma_{1e}}{L \mu}
		\begin{bmatrix}
			2\partial_{\overline{x}} \overline{u} & \frac{1}{\varepsilon}\partial_{\overline{y}} \overline{u} + \varepsilon \partial_{\overline{x}} \overline{v}\\
			\frac{1}{\varepsilon}\partial_{\overline{y}} \overline{u} + \varepsilon\partial_{\overline{x}} \overline{v} & 2 \partial_{\overline{y}} \overline{v}
		\end{bmatrix}.
	\end{equation*}
	Therefore, the normal and the tangential stress on the surface $\overline{y}=\overline{h}(\overline{x}, \overline{t})$ are given by,
	\begin{equation}
		\label{normal_stress}
		\renewcommand{\arraystretch}{1.5}
		\left[\mathbb{D}\vc{v}^L\cdot \vc{n}_1\right]\cdot \vc{n}_1 = \frac{\varepsilon^3 \sigma_{1e}}{L \mu}\frac{1}{\left( 1+ \varepsilon^2|\partial_{\overline{x}} \overline{h}|^2\right) } \left( - \partial_{\overline{y}} \overline{u} \;\partial_{\overline{x}} \overline{h} + \partial_{\overline{y}} \overline{v} + O(\varepsilon^2)\right)
	\end{equation}
	and
	\begin{equation}
		\label{tangential_stress}
		\renewcommand{\arraystretch}{1.5}
		\left[ \mathbb{D}\vc{v}^L\cdot \vc{n}_1\right]\cdot \vc{\tau} = \frac{\varepsilon^3\sigma_{1e}}{L \mu}\frac{1}{2}\frac{1}{\left( 1+ \varepsilon^2|\partial_{\overline{x}} \overline{h}|^2\right) } \left(- \frac{1}{\varepsilon} \partial_{\overline{y}} \overline{u} + O(\varepsilon)\right) .
	\end{equation}
	Similarly, the tangential gradient on $\overline{y}=\overline{h}(\overline{x}, \overline{t})$ becomes (cf. (\ref{tan_grad})),
	\begin{equation}
		\label{tan_grad1}
		\nabla_{\tau} p^s_1  
		= \frac{1}{L \;(1+ \varepsilon^2|\partial_{\overline{x}} \overline{h}|^2)^{1/2} } \partial_{\overline{x}} \tilde{p}^s_1.
	\end{equation}
	For the differential terms defined along the surface such as $\partial_t^\Sigma, \mathrm{div}_\Sigma$, we can compute them using the explicit forms in Appendix \ref{explicit_terms}. Also we will not distinguish below between $\tilde{p}^s_1$ and $p^s_1$, if not necessary.
	
	On the liquid-solid surface $\overline{y}=0$, corresponding quantities are simpler. The outward unit normal vector with respect to the solid surface is $\vc{n}_2=(0,1),$ whereas $\tau = (1,0)$, $\kappa_\Sigma =0$. Thus, $\left[ \mathbb{D}\vc{v}^L\cdot \vc{n}_2\right]\cdot \vc{\tau} = \frac{1}{2}\frac{\varepsilon^3 \sigma_{1e}}{L\mu}\left( \frac{1}{\varepsilon}\partial_{\overline{y}} \overline{u} + \varepsilon \partial_{\overline{x}} \overline{v}\right) $, and $\nabla_{\tau}p^s_2 = \frac{1}{L} \; \partial_{\overline{x}} \tilde{p}^s_2$.
	
	We also use the linear approximation for the surface pressures (\ref{8}) and (\ref{16}) in the following form involving the equilibrium density functions $\varrho^s_{ie}$,
	\begin{equation}
		\label{surface_tension}
		\sigma^s_i\equiv -p^s_i =\sigma_{ie}+ \gamma_i (\varrho^s_{ie}-\varrho^s_i) , \qquad i=1,2,
	\end{equation}
	where $ \sigma_{ie} := \sigma^s_i ( \varrho^s_{ie}) = -p^s_i ( \varrho^s_{ie}) = - \gamma_i (\varrho^s_{ie} - \varrho^s_{i,0})>0$.
	
	Further, the following assumptions are made,
	\begin{equation}
		\label{assump}
		\tau_1, \tau_2 \ll \frac{L \mu}{\varepsilon^3\sigma_{1e}} \quad \text{ and } \quad \frac{\varepsilon^2\sigma_{1e}}{\gamma_1} \ll \varrho^s_{1e}, \varrho^s_{2e},
	\end{equation}
	where the first condition means that the relaxation times $\tau_i$ are very small compared to the macroscopic time scale and the second condition says the deviation of $\varrho^s_i$ from its equilibrium value $\varrho^s_{ie}$, on both interfaces, is also small.
	
	\paragraph{Fluid equation:}
	Under the scaling, the incompressibility condition $(\ref{1})_1$ remains same in the new variables, leading to $(\ref{20})_1$. The Navier-Stokes equation $(\ref{1})_2$ becomes in the first component, with $\mathrm{Re}= \mu^2/\varrho^L \varepsilon^3 \sigma_{1e} L$
	\begin{equation*}
		\mathrm{Re} \left( \partial_{\overline{t}} \overline{u} + \overline{u} \partial_{\overline{x}} \overline{u} + \overline{v} \partial_{\overline{y}} \overline{u}\right) =  \Big( \partial_{\overline{x}}^2\overline{u} + \frac{1}{\varepsilon^2} \partial_{\overline{y}}^2\overline{u}\Big) - \frac{1}{\varepsilon^2} \partial_{\overline{x}} \overline{p}.
	\end{equation*}
	Thus in the limit $\varepsilon \to 0$, both the time derivative and the non-linear term vanish compared to the viscous term when Reynolds number satisfies $\varepsilon^2 \textrm{Re}
	\ll 1$ and one obtains $(\ref{20})_2$.
	Similarly the second component of the Navier-Stokes equation reduces to $(\ref{20})_3$.
	
	\paragraph{Boundary conditions:}
	The kinematic condition (\ref{2}) becomes in the limit (\ref{21.1}), using the fact that the new free boundary is $\overline{y}=\overline{h}(\overline{x},\overline{t})$ (cf. (\ref{normal_1})). 
	For the other boundary conditions, writing $\varrho^s_i = \varrho^s_{ie} + \frac{\varepsilon^2 \sigma_{1e}}{\gamma_1} \zeta_i, i=1,2$, one can compute the terms of equations (\ref{3}) and (\ref{9}) (cf. Appendix \ref{explicit_terms}) which yields,
	\begin{equation*}
		\frac{\varepsilon^3 \sigma_{1e}}{L \mu} \; \frac{\varepsilon^2 \sigma_{1e}}{\gamma_1} \;\; \partial_{\overline{t}} \zeta_i + \frac{1}{L} \frac{\varepsilon^3 \sigma_{1e}}{\mu} \partial_{\overline{x}} \Big( \overline{u}^s_{i}\Big( \varrho^s_{ie} + \frac{\varepsilon^2 \sigma_{1e}}{\gamma_1} \zeta_i\Big) \Big) = - \frac{1}{\tau_i} \frac{\varepsilon^2 \sigma_{1e}}{\gamma_1} \zeta_i.
	\end{equation*}
	Taking into account assumptions (\ref{assump}), we obtain the limit equations (\ref{21.2}) and (\ref{23.1}),
	where the time derivative vanishes and we assume $ \frac{1}{\lambda_i} :=\frac{L\mu}{ \varepsilon \tau_i \gamma_1 \varrho^s_{ie}} = O(1)$.
	Equation (\ref{4}) (similarly (\ref{12})) reduces to,
	\begin{equation}
		\label{10}
		\left( \overline{u}-\overline{u}^s_{1} \right) \partial_{\overline{x}} \overline{h} =\left( \overline{v}- \overline{v}^s_{1}\right) + \frac{\mu}{\varepsilon^2 \gamma_1 \tau_1 \varrho^L}\zeta_1.
	\end{equation}
	At this point, note that, as the thickness of the interface $\delta$ is assumed to be very small compared to the height of the fluid, i.e. $\delta \ll \varepsilon L$, and we also assume that the surface density $\varrho^s_i$ is of the same order as the bulk fluid density $\varrho^L$, i.e. $\varrho^s_i L^2 \sim \delta L^2 \varrho^L$, we have,
	$
	\varrho^s_i \ll \varepsilon L \varrho^L.
	$
	With the assumption that $ \frac{L\mu}{ \varepsilon \tau_1 \gamma_1 \varrho^s_{1e}} = O(1)$, the last term of equation (\ref{10}) then vanishes in the limit and we obtain (\ref{21.3}) (similarly (\ref{23.2})).
	To rewrite equations (\ref{5})-(\ref{7}), (\ref{13}), (\ref{15}), one needs to use approximation (\ref{8}) and (\ref{16}), or equivalently (\ref{surface_tension}) which says
	\begin{equation}
		\label{surface_tension1}
		p^s_1 = \sigma_{1e} \left( \varepsilon^2 \zeta_1 -1\right) \qquad \text{ and } \quad p^s_2 = \varepsilon^2 \sigma_{1e}\frac{\gamma_2}{ \gamma_1} \zeta_2 - \sigma_{2e}.
	\end{equation}
	Then equation (\ref{5}) becomes, with the help of (\ref{normal_stress}), (\ref{surface_tension1}) and (\ref{curvature}),
	\begin{equation*}
		\begin{aligned}
			&p^G - \frac{\varepsilon \sigma_{1e}}{L} \overline{p} + 2 \frac{\varepsilon^3 \sigma_{1e}}{L}\frac{1}{\left( 1+ \varepsilon^2|\partial_{\overline{x}} \overline{h}|^2\right) } \left( - \partial_{\overline{y}} \overline{u} \;\partial_{\overline{x}} \overline{h} + \partial_{\overline{y}} \overline{v} + O(\varepsilon^2)\right)\\
			& \hspace{4cm}+ \sigma_{1e} \left( \varepsilon^2 \zeta_1 -1\right) \frac{\varepsilon \;\partial^2_{\overline{x}} \overline{h}}{L\left( 1+ \varepsilon^2 |\partial_{\overline{x}} \overline{h}|^2\right) ^{3/2}} =0,
		\end{aligned}
	\end{equation*}
	which one can write in the following form,
	\begin{equation*}
		\begin{aligned}
			& L \varepsilon^2 p^G - \varepsilon^3 \sigma_{1e} \overline{p} + 2 \varepsilon^5 \sigma_{1e}\frac{1}{\left( 1+ \varepsilon^2|\partial_{\overline{x}} \overline{h}|^2\right) } \left( - \partial_{\overline{y}} \overline{u} \;\partial_{\overline{x}} \overline{h} + \partial_{\overline{y}} \overline{v} + O(\varepsilon^2)\right)\\
			& \hspace{4cm} + \left( \varepsilon^5 \sigma_{1e} \zeta_1 -\varepsilon^3 \sigma_{1e}\right) \frac{ \;\partial^2_{\overline{x}} \overline{h}}{\left( 1+ \varepsilon^2 |\partial_{\overline{x}} \overline{h}|^2\right) ^{3/2}} =0.
		\end{aligned}
	\end{equation*}
	Now, assuming that $\varepsilon^3 \sigma_{1e} = O(1)$, one obtains in the limit (\ref{22.1}).
	For equation (\ref{6}) (similarly (\ref{13})), plugging in the expressions (\ref{tangential_stress}), (\ref{tan_grad1}) and (\ref{surface_tension1}), one gets
	\begin{equation*}
		\frac{\varepsilon^3\sigma_{1e}}{L}\frac{1}{\left( 1+ \varepsilon^2|\partial_{\overline{x}} \overline{h}|^2\right) } \Big(- \frac{1}{\varepsilon} \partial_{\overline{y}} \overline{u} + O(\varepsilon)\Big) = \frac{\varepsilon^2\sigma_{1e}}{L \;(1+ \varepsilon^2|\partial_{\overline{x}} \overline{h}|^2)^{1/2} } \partial_{\overline{x}} \zeta_1,
	\end{equation*}
	which, taking into account the scaling $\varepsilon^3 \sigma_{1e} = O(1)$, simplifies to (\ref{22.2})
	to leading orders. Also, equation (\ref{7}) (similarly (\ref{15})) for the tangential surface velocity becomes in the new variables, combining (\ref{tan_grad1}) and (\ref{surface_tension1}),
	\begin{equation*}
		-\frac{(1+4 \alpha_1 \beta_1)}{4 \beta_1} \frac{\varepsilon^2\sigma_{1e}}{L \;(1+ \varepsilon^2|\partial_{\overline{x}} \overline{h}|^2)^{1/2} } \partial_{\overline{x}} \zeta_1 = \frac{\varepsilon^3 \sigma_{1e}}{\mu} \frac{\left( \left( \overline{u}^s_1 - \overline{u}\right) + \varepsilon^2 \partial_{\overline{x}} \overline{h} \left( \overline{v}^s_1 - \overline{v}\right)\right)}{(1+ \varepsilon^2 |\partial_{\overline{x}} \overline{h}|^2)^{1/2}},
	\end{equation*}
	which reduces to (\ref{22.3}),
	where we make the assumption that $\frac{\mu}{L}\frac{(1+4 \alpha_1 \beta_1)}{4 \varepsilon \beta_1} = O(1)$.
	Lastly, equation (\ref{14}) can be computed easily in the same manner.
	
	\paragraph{Contact line condition:}
	For the equations at the contact line, the velocity of the contact point can be calculated as $\vc{v}^c = \partial_t s_i= \frac{\varepsilon^3 \sigma_{1e}}{\mu}\left( \partial_{\overline{t}}\overline{\Lambda}_i(\overline{t}),0\right), i=1,2$. Also, the contact angle being very small, the normal vectors at the contact point satisfy $\nu_1 = \nu_2= (1,0)$. Therefore, the continuity equation (\ref{17}) reads as
	\begin{equation*}
		\Big( \varrho^s_{1e} + \frac{\varepsilon^2 \sigma_{1e}}{\gamma_1} \zeta_1\Big) \frac{\varepsilon^3 \sigma_{1e}}{\mu}\left(\overline{u}^s_1 -\partial_{\overline{t}}\overline{\Lambda}_1\right)  + \Big( \varrho^s_{2e} + \frac{\varepsilon^2 \sigma_{1e}}{\gamma_1} \zeta_2\Big)\frac{\varepsilon^3 \sigma_{1e}}{\mu}\left(\overline{u}^s_2 -\partial_{\overline{t}}\overline{\Lambda}_2\right) =0,
	\end{equation*}
	which simplifies to (\ref{44.1}), taking into account assumption (\ref{assump}).
	Also the contact angle being small, we approximate $\cos\theta \approx 1- \frac{\theta^2}{2} $. Further, let us assume that
	$p^s_3 = \sigma_{3e}\left( 1- \varepsilon^2 \overline{\sigma}_3\right) \ge 0$ is a constant. Note that the surface tension on the solid-gas interface works in opposite direction to that of the liquid-solid surface. Therefore, as the surface tensions at equilibrium are balanced, i.e. \begin{equation}
		\label{equi_surface_tension}
		\sigma_{1e} + \sigma_{2e}+ \sigma_{3e} =0,
	\end{equation}
	the Young's equation (\ref{18}) can be written as, using (\ref{surface_tension1}),
	\begin{equation*}
		\sigma_{1e} \left( \varepsilon^2 \zeta_1 -1\right) \Big( 1- \frac{\varepsilon^2}{2} |\partial_{\overline{x}} \overline{h}|^2\Big) = \sigma_{3e}\left( 1- \varepsilon^2 \overline{\sigma}_3\right) - \Big( \varepsilon^2 \sigma_{1e}\frac{\gamma_2}{ \gamma_1} \zeta_2 - \sigma_{2e}\Big).
	\end{equation*} 
	Using the equilibrium condition (\ref{equi_surface_tension}), the above contact angle condition simplifies to (\ref{44.3}).
	Finally, the missing condition (\ref{18.1}) becomes (\ref{44.4}) with standard computation. From now on, we omit the bars for convenience.
	
	\subsubsection{Derivation of thin film equation}
	
	Next, we would like to derive the lubrication equation for the profile $h(x,t)$ from equations (\ref{20})-(\ref{44.4}). We show below that with $g = \frac{\varrho^s_{2e} \tau_2 }{\varrho^s_{1e} \tau_1}, b_1= a\left( \frac{1}{4} + \alpha_2 \beta_2\right),  \frac{\varrho^s_{1e}}{\varrho^s_{2e}},$ $c_1 = \left( 1-2\beta d_1 q\right),$ $c_2 = 2 a\left( \frac{1}{4}-\alpha_2 \beta_2\right)$, for $t>0$, one obtains,\\
	\\
	\begin{equation}
		\label{thin_film_model}
		\left.
		\begin{aligned}
			\partial_{t} h +\partial_{x}\Big( \Big( \frac{h}{3} + \frac{1}{\beta}\Big) h^2\;\partial^3_{x} h - \frac{1}{2} \partial_{x} \zeta_1 \;h^2 -\frac{1}{\beta} \Big(\partial_{x} \zeta_1 + \frac{a}{2} \;\partial_{x} \zeta_2\Big) h\Big) &=0 ,\\
			\lambda_1\;\partial_{x} \Big((h+ \frac{1}{\beta}+ d_1) \partial_x \zeta_1 + \frac{a}{2 \beta}\;\partial_x \zeta_2 - \Big( \frac{h}{2}+ \frac{1}{\beta}\Big) h\;\partial_x^3 h\Big) &=\zeta_1,\\
			\frac{g \lambda_1}{ \beta} \;\partial_{x} \Big(b_1 \;\partial_x \zeta_2 +\frac{1}{2} \left( \partial_x \zeta_1 - h\;\partial^3_{x} h \right)\Big) &=\zeta_2,
		\end{aligned}
		\right\}
		x\in (\Lambda_1, \Lambda_2);
	\end{equation}
	
	\begin{equation}
		\label{thin_film_bc}
		\left.
		\begin{aligned}
			- h \;\partial_x^3 h + c_1 \; \partial_x \zeta_1 + c_2 \; \partial_x \zeta_2 &= 0,\\
			\zeta_1 +a\; \zeta_2 +\frac{1}{2} |\partial_{x} h|^2 &=b,\\
			m \;\zeta_1 &=\zeta_2,\\
			h &=0 .
		\end{aligned}
		\right\}
		\quad x= \{\Lambda_1, \Lambda_2\}.
	\end{equation}
	
	Indeed, from $(\ref{20})_3$ and (\ref{22.1}), we can first conclude that $p$ is independent of $y$. 
	Thus, relation $(\ref{20})_2$ gives
	\begin{equation}
		\label{32}
		u(x,y) = \frac{1}{2} \partial_{x} p
		\left( y - 2 h\right) {y}+ A(x) y + B(x) \qquad \text{ in } \ (\Lambda_1, \Lambda_2)\times (0,h),
	\end{equation}
	which yields, due to relations (\ref{22.2}) and (\ref{31.1}) that
	\begin{equation*}
		A(x) = -\partial_{x} \zeta_1 \quad \text{ and } \quad B(x) = -\frac{1}{ \beta}\left( \partial_{x} \zeta_1 + \partial_{x} p\;h + \frac{a}{2} \partial_{x} \zeta_2\right) .
	\end{equation*}
	Thus, (\ref{32}) becomes, together with (\ref{22.1}),
	\begin{equation}
		\label{u}
		u(x,y) = \frac{1}{2 } \partial^3_{x} h\left( 2h-y\right) y -\partial_{x} \zeta_1 \;y -\frac{1}{ \beta}\left( \partial_{x} \zeta_1 - h\;\partial^3_{x} h + \frac{a}{2} \partial_{x} \zeta_2\right) \quad \text{ in } (\Lambda_1, \Lambda_2)\times (0,h).
	\end{equation}
	Now, since due to (\ref{31.3}) and (\ref{23.2}), we have $v\arrowvert_{y= 0}=0$, one gets from the divergence condition $(\ref{20})_1$,
	\begin{equation*}
		v\arrowvert_{y= h} = -\int_0^{h}{ \partial_{x} u \;\mathrm{d}y} = - \partial_{x} \left( \int_0^{h}{u \;\mathrm{d}y}\right) + u\;\partial_{x} h.
	\end{equation*}
	Combining this relation with (\ref{21.1}) and (\ref{21.3}), we can write the equation for the normal velocity of the free surface in the following conservation form,
	\begin{equation}
		\label{24.}
		\partial_{t} h= - \partial_{x} \Big( \int_0^{h}{u \;\mathrm{d}y}\Big) .
	\end{equation}
	Therefore, with the help of (\ref{u}), we obtain the thin film equation $(\ref{thin_film_model})_1$.
	Furthermore, plugging the value of $u^s_1$ from (\ref{22.3}) into the equation (\ref{21.2}), one obtains $(\ref{thin_film_model})_2$, using (\ref{u}) at $y=h$.
	Similarly, equation (\ref{23.1}) yields  $(\ref{thin_film_model})_3$, with the help of (\ref{31.2}) and (\ref{u}), as
	\begin{equation*}
		\begin{aligned}
			-\zeta_2 & = g \lambda_1\;\partial_{x} \Big( -\frac{1}{2 \beta}\Big(\partial_{x} \zeta_1 - \partial^3_{x} h
			\; {h} + \frac{a}{2}\partial_{x} \zeta_2 \Big) -d_2\; \partial_{x} \zeta_2 \Big)\\
			& = -\frac{g \lambda_1}{ \beta}\partial_{x} \Big( \frac{1}{2}\Big(\partial_{x} \zeta_1 - \partial^3_{x} h
			\; {h} \Big) + b_1 \; \partial_{x} \zeta_2 \Big) \qquad \text{ on } y=0.
		\end{aligned}
	\end{equation*}
	Finally, the mass balance condition at the triple junction (\ref{44.1}) gives, with the help of the equations (\ref{22.3}) and (\ref{31.2}) for $u^s_i$, together with the expression (\ref{u}) for $u$,
	\begin{equation}
		\label{26}
		\begin{aligned}
			&	-(1+q) \beta \; \partial_{t}\Lambda_i + q\Big[\partial^3_{x} h\; h - (1+\beta d_1)\partial_{x} \zeta_1 - \frac{a}{2} \;\partial_{x} \zeta_2 \Big] \\
			& \hspace{2.5cm} = \Big[ \frac{1}{2}\left( \partial_{x} \zeta_1-\partial^3_{x} h
			\; {h}\right) + \left( \frac{a}{4}+ \beta d_2\right) \partial_{x} \zeta_2 \Big] \qquad  \text{ at } y = h=0.
		\end{aligned}
	\end{equation}
	Note that for the equilibrium surface density in the liquid-solid surface, it holds $\varrho^s_{2e}> \varrho^s_{1e}$ which means $q<1$.
	Also noting that the velocity of the contact point $s_i(t)$ is given by, from the thin film equation $(\ref{thin_film_model})_1$,
	$$
	\partial_{t} \Lambda_i = \frac{1}{ \beta} \Big( h \;\partial^3_{x} h- \partial_{x} \zeta_1 - \frac{a}{2}\; \partial_{x}\zeta_2 \Big) \Big\arrowvert _{x=\Lambda_i},
	$$
	we may further simplify the above boundary condition (\ref{26}) to $(\ref{thin_film_bc})_1$.
	
	Observe that in the final thin film model (\ref{thin_film_model})-(\ref{thin_film_bc}) we have only two parameters: $\frac{1}{\beta}$ which represents the ratio of the slip-length and the film thickness and $\lambda_1$ which is of order $1$ or small.
	
	\begin{remark} 
		In the case when the slip length is small compared to the film height, i.e. $\frac{1}{\beta}\ll 1$, we obtain from $(\ref{thin_film_model})_1$,
		\begin{equation*}
			\partial_{t} h +\partial_{x}\Big( \tfrac{1}{3} h^3\;\partial^3_{x} h - \tfrac{1}{2} \partial_{x} \zeta_1 \;h^2 \Big) =0 ,
		\end{equation*}
		which can be seen as the counterpart of the classical thin film equation with no-slip boundary condition
		\begin{equation*}
			\partial_t h + \partial_x\big( h^3\partial^3_x h\big) =0.
		\end{equation*}
		
	\end{remark}
	
	\subsection{Thin film model for large slip length}
	\label{S2.2}
	
	If we assume that the slip length is large compared to the film thickness i.e. $\frac{1}{\beta}\gg 1$, one obtains from (\ref{thin_film_model})-(\ref{thin_film_bc}), with the new time scale $\tilde{t} 
	= \frac{ \varepsilon^3 \sigma_{1e}}{L \mu\beta}t$ and the variable $\tilde{\lambda}_1 =\frac{\lambda_1}{\beta}=O(1)$, the following thin film model, for $t>0$, (removing the tildas once again)
	\begin{equation*}
		\left.
		\begin{aligned}
			\partial_{t} h +\partial_{x}\Big( h^2\partial^3_{x} h -\Big(\partial_{x} \zeta_1 + \tfrac{a}{2} \;\partial_{x} \zeta_2\Big) h\Big) &=0,\\
			\lambda_1\;\partial_{x} \Big( \partial_x \zeta_1 + \tfrac{a}{2} \;\partial_x \zeta_2 - h\;\partial_x^3 h\Big) &=\zeta_1,\\
			g \lambda_1 \;\partial_{x} \Big(\tfrac{a}{4}\partial_x \zeta_2 +\tfrac{1}{2} \Big( \partial_x \zeta_1 - h\;\partial^3_{x} h \Big)\Big) &=\zeta_2,
		\end{aligned}
		\right\}
		\quad	x\in (\Lambda_1, \Lambda_2);
	\end{equation*}
	\begin{equation*}
		\left.
		\begin{aligned}
			- h \;\partial_x^3 h + \partial_x \zeta_1 +\tfrac{a}{2} \; \partial_x \zeta_2 &= 0,\\
			\zeta_1 + a\;\zeta_2 +\tfrac{1}{2} |\partial_{x} h|^2 &= b,\\
			m\;\zeta_1 &=\zeta_2,\\
			h &=0,
		\end{aligned}
		\right\}
		\quad	x= \{\Lambda_1, \Lambda_2\}.
	\end{equation*}
	
	\begin{remark}
		The above thin film equation can be compared to the standard one with slip boundary condition,
		\begin{equation*}
			\partial_t h + \partial_x\left( h^2\partial^3_x h\right) =0.
		\end{equation*}
	\end{remark} 
	
	\subsection{A solid wedge entering a thin liquid film}
	\label{S2.3}
	
\begin{figure}[h!]
		\centering
		\begin{minipage}{.33\textwidth}
			\centering
			\includegraphics[width=.6\linewidth]{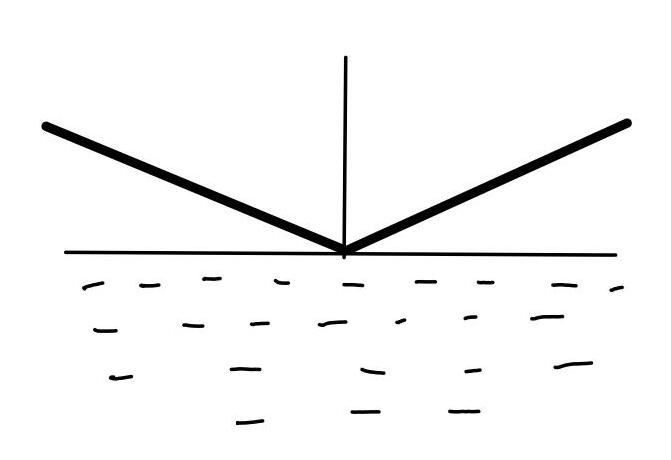}
		\end{minipage}%
		\begin{minipage}{.33\textwidth}
			\centering
			\includegraphics[width=.6\linewidth]{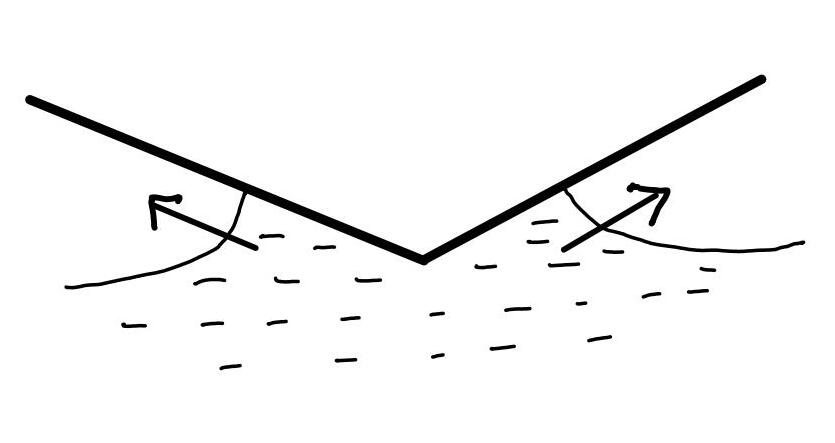}
		\end{minipage}
		\begin{minipage}{.33\textwidth}
			\centering
	\includegraphics[width=.6\linewidth]{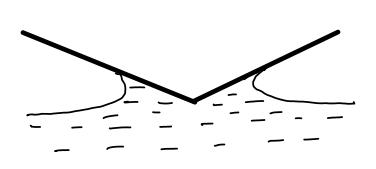}
		\end{minipage}
		\caption{General meniscus formation}
		\label{fig3}
	\end{figure}

	In this subsection, we briefly discuss a different thin film model, where a  wedge-shaped solid enters a thin liquid film  (cf. Figure \ref{fig3}) which leads to movement of  the contact line and the formation of  a meniscus. 
	We briefly discuss here the thin film model in this situation, both, in the standard approach, i.e. with slip-boundary conditions, and in the approach of Shikhmurzaev. We do not give here the details of the derivation, 
	 but rather state the limit models and discuss the main differences. 
	
	\FloatBarrier
	We consider a (horizontal, two-dimensional) thin fluid film of thickness $H$ and a wedge-shaped solid which initially touches the film with an angle $\delta$  (cf. Figure \ref{fig4}) and then enters the liquid such that
	$\tilde h(t)= g( \frac{t}{t_0}\big)^\eta $ with $\eta >0$, where $\tilde h(t)$ denotes the distance between the tip of the wedge and the bottom of the liquid domain. Here  $t_0$ is the time scale that characterises the motion of the solid wedge moving down. The particular case of constant velocity of the wedge corresponds to $\eta=1$, while the case of constant acceleration corresponds to $\eta=2$.
	As a consequence, the contact point moves and makes an angle $(\pi -\theta)$ between the liquid-solid and the liquid-gas interface. The liquid-gas interface is described by $y=h(x,t)$ and the liquid-solid rigid boundary is given by $y=\tilde{h}(t) +\delta |x| $. The free interface makes an angle $\tilde{\theta}$ with the horizontal plane (cf. Figure \ref{fig5}). Therefore, $\theta = \tilde{\theta}+ \delta$ and $\tan\tilde{\theta} = -\partial_x h$. To obtain the lubrication approximation, we need to assume that the two angles $\theta$ and $\delta$ are very small and of the same order, which means that the contact angle $(\pi -\theta)$ is close to $\pi$ in this setting. Note that the thin film approximation would not be valid if $\delta$ is large. The contact points are given by $s_i = (\Lambda_i(t), y)$.
	
	Finally, we assume that the far field condition $h \to H$ as $|x|\to \infty$ holds.
	\begin{figure}[h]
		\centering
		\begin{minipage}{.45\textwidth}
			\centering
			\begin{tikzpicture}
				\draw[-,thick] (-2,0)--(2,0);
				\draw[-,thick] (-2,-1)--(2,-1);
				\draw[-,thick] (0,0)--(2,1);
				\draw[-,thick] (0,0)--(-2,1);
				\draw[<->] (1.9,-0.05)--(1.9,-0.95);
				\node at (2.05, -0.5) {$\varepsilon$};
				\draw[->] (-1.9,0.7)--(-1.6,0.7);
				\draw (0.49,0) arc (0:60:0.22);
				\node at (-2.37,0.7) {solid};
				\node at (2, 0.45) {gas};
				\node at (0,-0.6) {liquid};
				\node at (0.75, 0.18) {$\delta$};
			\end{tikzpicture}
			\subcaption{} \label{fig4}
		\end{minipage}
		\begin{minipage}{.49\textwidth}
			\centering
			\begin{tikzpicture}
				\draw[dashed] (-2.5,0)--(2.5,0);
				\draw[-,thick] (-2.5,-1)--(2.5,-1);
				\draw[-,very thick] (0,-0.3)--(1,0.22);
				\draw[dashed] (0.5,0.22)--(1.7,0.22);
				\draw[-,very thick] (0,-0.3)--(-1,0.22);
				\draw (0.57,0)--(2.5,1);
				\draw (-0.57,0)--(-2.5,1);
				\draw[<->] (2.45,-0.05)--(2.45,-0.95);
				\node at (2.58, -0.5) {$\varepsilon$};
				\draw[<->] (0,-0.3) -- (0,-0.95);
				\node at (-0.2, -0.6) {$\tilde{h}$};
				\draw[dashed] (1,0.21) -- (1,-0.95);
				\node at (1, -1.3) {$\Lambda$};
				\draw[dashed] (-1,0.21) -- (-1,-0.95);
				\draw[thick] (1,0.22) .. controls (1.4,0) and (2,0) .. (2.5,0);
				\draw[<-] (2.1,0.05)--(2.4,0.26);
				\node at (3.4,0.35) {$y=h(x,t)$};
				\draw[thick] (-1,0.22) .. controls (-1.4,0) and (-2,0) .. (-2.5,0);
				\draw (1.45,0.09) arc (0:50:0.16);
				\draw[<-] (1.55,0.15)--(1.85,0.35);
				\node at (1.98, 0.45) {$\tilde{\theta}$};
			\end{tikzpicture}
			\subcaption{}
			\label{fig5}
		\end{minipage}
		\caption{}
	\end{figure}
	
	Since we are in a symmetric configuration with respect to the $y$-axis, it is enough to analyse the situation $x>0$. 
	
	\subsubsection{Classical approach}
	\label{S2.3.1}
	
	In classical fluid mechanics, conditions at contact lines are considered only at equilibrium, hence the contact angle is prescribed. Starting from the incompressible Navier-Stokes equation
	in the liquid with a slip boundary condition at the fluid-solid interface, one can derive, after an appropriate rescaling, the following thin film equation 
	\begin{align}
		\partial_{t} h + \partial_{x}\Big(  \Big( \frac{h}{3} + \frac{1}{\beta} \Big)h^2 \partial^3_{x} h \Big)  =0, \qquad &\text{ in } x> \Lambda, \label{thin_film_eqn}\\
		h= \Big( 1 - \Big( \frac{t}{t_0}\Big)^\eta \Big) + k_1 \Lambda, \quad
		\partial_{x} h = k_1 -k_2, \qquad &\text{ at } x = \Lambda, \label{contact_line} \\
		h =1, \qquad &\text{ at } \ t=0, \label{IC}
	\end{align}
	where $k_1,$ and $k_2$ are rescalings of the angles $\delta$ and $\theta$, 
	together with the relation
	\begin{equation}
		\label{extra_cond}
		h \;\partial_{x}^3 h = A \qquad \text{ at } x = \Lambda.
	\end{equation}
	where  $A$  can be computed from the inner flow via solving a linear ODE.
	
	Relation (\ref{extra_cond}), together with the thin film equations (\ref{thin_film_eqn})-(\ref{IC}) determines fully the dynamics and the position of the contact point.
	
	\subsubsection{Shikhmurzaev's approach}
	\label{S2.3.2}
	
	For  Shikhmurzaev's approach, 
	we use the same scaling as in Section \ref{S2.1}, only with the adaptation for the quantities $\tilde{h}, t_0,$ and $\delta$, that is,
	\begin{equation*}
		\overline{\tilde{h}} = \frac{\tilde{h}}{H}, \qquad \overline{t}_0 = \frac{\varepsilon^3 \sigma_{1e}}{L \mu} t_0, \qquad \delta = k_1 \varepsilon.
	\end{equation*}
	
	The difference with the thin film model described in Section \ref{S2.1} or \ref{S2.2} is that, the equations on the liquid-gas interface (\ref{2})-(\ref{8}) are valid for $y=h(x,t), x> \Lambda(t)$, while the liquid-solid interface conditions (\ref{9})-(\ref{16}) are valid for $y= \tilde{h}(t) + \delta x, x\in (0,\Lambda(t))$, with $\vc{v}^S = (0, \partial_t \tilde{h})$ (at the upper boundary) and for $y=0$ with $\vc{v}^S=0$ (at the lower boundary). Using the same scaling as in Section \ref{S2.1}, we only indicate here the differences from the previous classical case. Essentially, it is enough to calculate the boundary conditions at the liquid-solid interface $y= \tilde{h}(t) + \delta x, x\in (0,\Lambda(t))$ together with the contact line conditions (\ref{17})-(\ref{18.1}).
	
	We obtain as before, equations (\ref{thin_film_model})
	in $t>0, \; x>\Lambda$.
	On $ t>0$ and $x= \Lambda$, we obtain the following boundary conditions (cf. (\ref{thin_film_bc}) where $(\ref{thin_film_bc})_4$ is replaced by $(\ref{contact_line})_1$),
	\begin{equation}
		\label{thin_film_bc_wedge}
		\begin{aligned}
			-\frac{1}{2} h \;\partial_x^3 h + c_1\; \partial_x \zeta_1 + c_2 \; \partial_x \zeta_2 &= 0,\\
			\zeta_1 +a\; \zeta_2 +\frac{1}{2} |\partial_{x} h|^2 &=b,\\
			m \;\zeta_1 &=\zeta_2,\\
			h &= \Big( 1 - \Big( \frac{t}{t_0}\Big)^\eta \Big) + k_1 \Lambda
		\end{aligned}
	\end{equation}
	together with the equation (due to the continuity of velocity $u$ at the contact point)
	\begin{equation}
		\label{48}
		A = h \;\partial_x^3 h - \partial_x \zeta_1 \qquad \text{ at } x = \Lambda,
	\end{equation}
	Here $A$ is determined as follows: it solves,
	\begin{equation}
		\label{A1}
		\tilde{r}_1(x, t) \partial_x A + \tilde{r}_2(x, t) A + F\left( \partial_x^2 \zeta_2, \partial_x \zeta_2\right) + \partial_t \tilde{h}=0, \qquad x \in (0, \Lambda),
	\end{equation}
	where
	\begin{equation*}
		\tilde{r}_1(x,t) = ( \tilde{h}+ k_1 x) \Big( \frac{1}{\beta} + \frac{1}{6}( \tilde{h}+ k_1 x )\Big), \quad \tilde{r}_2(x,t) = k_1 \Big( \frac{1}{\beta} + \frac{1}{3}( \tilde{h}+ k_1 x)\Big) ,
	\end{equation*}
	and
	\begin{equation*}
		F = -\frac{a}{2\beta} \big( k_1 \partial_x \zeta_2 + ( \tilde{h}+ k_1 x) \partial^2_x \zeta_2\big) .
	\end{equation*} 
	Also the surface density $\zeta_2$ satisfies, in the inner region,
	\begin{equation}
		\label{zeta}
		-\zeta_2 = \lambda_2\; \partial_x\Big[-\frac{1}{4} ( \tilde{h}+ k_1 x)^2 
		\frac{2A}{( \tilde h + k_1 x)} + \Big( \tilde{h}+ k_1 x + \frac{1}{ \beta}\Big) A - \Big( \frac{a}{2\beta}+ d_2\Big) \partial_x \zeta_2 \Big], \quad x\in(0,\Lambda) 
	\end{equation}
	The symmetry assumption of the configuration further gives $u=0$ at $x=0$ and thus,
	\begin{equation}
		\label{BC}
		A =0, \quad \partial_x \zeta_2=0 \qquad \text{ at } x =0.
	\end{equation}
	Equations (\ref{A1}) and (\ref{zeta}) are a coupled system with boundary conditions (\ref{BC}). Solving this would provide $A$ at $x=\Lambda$ which in turn gives us the necessary condition (\ref{48}) in order to determine the position of the contact point.
	
	Summing up, we obtain the thin film equation (\ref{thin_film_model}), (\ref{thin_film_bc_wedge}) together with the condition (\ref{48}).
	
\section{Conclusion and discussion}
	
In this paper, we rederive Shikhmurzaev's model for fluid interfaces. An appealing aspect of the model is that it provides a set of equations describing the variables that characterizes the state of the interfaces. In particular, we study in detail under which conditions it is possible to derive isothermal limits. We also conclude, as it was noticed in some previous works (\cite{bil06}, \cite{bothe}) that there is a missing boundary condition at the triple junction in the original formulation of Shikhmurzaev model. We further derive in this paper some thin film approximations that provide simpler PDE models than the original one.


	
\paragraph{Physical assumptions of Shikhmurzaev's model}


One of the interesting questions is to understand if Shikhmurzaev's model makes any
kind of experimental prediction that differs from the traditional
models. It is difficult to get a clear insight
into the model with so many parameters, even in selected asymptotic regimes. The real values of different parameters are also hard to obtain and compare. We think that the main physical feature that is
contained in Shikhmurzaev's model, but not in the classical
ones (including Navier slip boundary conditions, precursors, disjoining
pressure etc), is that it contains a characteristic
time scale that yields the time in which the interface moves
towards its equilibrium value. This time scale is assumed to be very short
($10^{-6}$s). Therefore, the effects associated to that scale would be
visible only in very fast phenomena. Shikhmurzaev describes the early stages
of
situations in which a solid object sinks in a fluid (see the discussion below). In the
very early stages of this process, the time scale introduced by Shikhmurzaev
might yield a different behaviour than the classical models. However, for
longer time scales, there might not be any difference between the
behaviour of these different models.

In this direction, a possible experiment that could illustrate the differences
between Shikhmurzaev models and classical models is the following (although we did not
workout the details further). The idea is to study the dynamics of oscillating bubbles, as they
appear in sonoluminiscence experiments. Bubbles in a sound field (with frequencies of at least 100,000 Hz) oscillate and the oscillating bubbles have a changing interface. At those
frequencies, the presence of the time scale for the relaxation of the
interfaces could interfere with the oscillations. We expect to have
different behaviours for the classical models and for Shikhmurzaev models in this setup.
Due to the spherical symmetry and the lack of contact
lines, it seems possible to analyze the model in detail. 
The analysis of this model would be particularly relevant in order to assert the validity of Shikhmurzaev's model in situations in which the dynamics of the bubble remains as simple as possible (for instance, region in which the spherical shape of the particle retains, shock waves do not arise etc).
	
\paragraph{Thin film model}
\label{conclusion_thin}
	
	From the discussion in Section \ref{S2.3.1} and Section \ref{S2.3.2} on the thin film model of a solid wedge moving downward into a fluid, with two different approaches, we may conclude the following. The classical model says that the contact angle is created instantaneously (cf. equation $(\ref{contact_line})_2$), while Shikhmurzaev's approach says that the contact angle is created in a continuous manner (cf. equation $(\ref{thin_film_bc_wedge})_2$).
	This above model of meniscus formation as a general phenomena (not the thin film approximation) has been mentioned in \cite[Fig 6]{shikh2020}.

	\renewcommand{\thetable}{\Alph{section}\arabic{table}}
	\appendix
	\section{Appendix}
	\label{appendix}
	
	\subsection{Notation}
	Throughout this work, we use bold letters to denote vectors (tensor of order one) and calligraphic letters to denote tensors (of second order).
	The exterior product of two vectors $\vc{u} = (u_l)_l, \vc{w}=(w_l)_l$  is defined as 
	$
	\vc{u} \otimes \vc{w} = \left( u_k w_l\right) _{kl}
	$
	while the  interior product of two tensors, indicated by a dot in between, is given with 
	$\mathbb{A}=(a_{kl})$ by
	$$
	\vc{u}\cdot \vc{w} = \sum_{l}u_l w_l , \quad 
	\mathbb{A}\cdot \vc{u} = \Big( \sum_{l} a_{kl} u_l\Big) _k, \quad \vc{u} \cdot \mathbb{A} = \Big (\sum_{l} u_l \;a_{lk}\Big) _k , \quad  
	\mathbb{A}\cdot \mathbb{B} = \Big( \sum_{l} a_{jl}b_{lk}\Big) _{jk}.
	$$
	In an analogous way, we define the scalar product of two tensors as,
	$
	\mathbb{A}:\mathbb{B} = \sum_{k,l} a_{kl} b_{kl}.
	$
	The operator of spatial derivation can be thought of as a vector $\frac{\partial }{\partial \vc{x}} \equiv \nabla$ with which the following operators are defined as,
	$$
	\nabla_{\vc{x}} a = \left( \frac{\partial a}{\partial x_l}\right)_{l}, \quad \nabla_{\vc{x}} \vc{u} = \left( \frac{\partial u_k}{\partial x_l}\right)_{kl} ,\quad
	\div \vc{u}\equiv \nabla_{\vc{x}} \cdot \vc{u} = \sum_{l} \frac{\partial u_l}{\partial x_l}, \quad \div \mathbb{A}= \Big( \sum_{l}\frac{\partial a_{kl}}{\partial x_l}\Big) _k.
	$$
	
	For every vector field $\vc{w}$ defined on an interface $\Sigma$ with the unit normal vector $\vc{n}$, we decompose it in the normal and tangential components by means of
	\begin{equation*}
		\vc{w} = \vc{w}_{\tau} + w_n \vc{n} \qquad \ \text{ where} \quad \vc{w}_\tau = (I - \vc{n}\otimes \vc{n}) \cdot \vc{w} \ \text{ and } \ w_n = \vc{w} \cdot \vc{n}.
	\end{equation*}
	Note that, according to the above notation, $w_n$ is a scalar quantity whereas $\vc{w}_\tau$ is a vector.
	
	
	\subsection{Differential geometry}
	We use $\partial_t^\Sigma$ to denote the normal derivative which is defined as follows: suppose a family of surfaces $\Sigma(t)$ is parametrized as $\Sigma(t)= \{\vc{x}: \vc{x} = \sigma (\vc{y}, t), \vc{y}\in D \subset \mathbb{R}^2\}$ where $\sigma$ has the property
	\begin{equation*}
		\partial_t \sigma(\vc{y},t) = v^s_n\; \vc{n},
	\end{equation*}
with $v^s_n$ being the normal velocity of the surface.	The existence of such a parametrization for every smoothly evolving family of surfaces $\Sigma(t)$ is easy to obtain. We then define, for a generic function $\varphi$, defined over $\Sigma$,
	\begin{equation}
		\label{normal_derivative.}
		\partial_t^\Sigma \varphi(\vc{x},t):= \frac{\mathrm{d}}{\mathrm{d}t}\varphi (\sigma(\vc{y},t),t).
	\end{equation}
	
	The surface divergence $\mathrm{div}_\Sigma$ is defined, for vector fields $\vc{u}$ defined on $\Sigma$ with values in the tangent space $T\Sigma$, i.e. $\vc{u}(\vc{x}) \in T\Sigma(\vc{x}),$ for all $\vc{x}\in\Sigma$, as
	\begin{equation*}
		\left\langle \mathrm{div}_{\Sigma}\; \vc{u}, \varphi\right\rangle_\Sigma := 
		- \int_\Sigma{\vc{u} \cdot \nabla_\tau \varphi}, \qquad \varphi\in C^\infty_c (\Sigma),
	\end{equation*}
	where $\nabla_\tau \varphi = (\mathbb{I} - \vc{n}\otimes \vc{n})\cdot \nabla_{\vc{x}}\varphi$ and $\varphi$ is extended by a constant along the normal in a small neighbourhood of $\Sigma$.
	
	In the similar way, the surface divergence	$\mathrm{div}_{\Sigma}$ for second order tensors is defined, for $\mathbb{T}(\vc{x}) \in\ \mathbb{R}^3 \otimes T_{\vc{x}}\Sigma$, as
	\begin{equation}
		\label{surface_divergence.}
		\left\langle \mathrm{div}_\Sigma \mathbb{T}, \varphi \right\rangle := 
			- \int_\Sigma{\mathbb{T}: \nabla_{\tau} \varphi}, \qquad \varphi \in C^\infty_c(\Sigma,\mathbb{R}^3).
		\end{equation}
		
		For a general smooth interface $\Sigma$ which is expressed as $\Sigma= \{\vc{x}: \vc{x} = \Gamma (\vc{y}, t), \vc{y}\in D \subset \mathbb{R}^2\}$
		where $\Gamma$ is a parametrization, the velocity $\vc{v}^s$ of the interface is given by
		$
		\vc{v}^s = \frac{\partial \Gamma}{\partial t}(\vc{y},t)
		$ at each point of the interface.
		If $V_n$ denotes the normal velocity of the interface, then
		\begin{equation}
			\label{46}
			V_n = \vc{v}^s\cdot \vc{n}.
		\end{equation}
		Note that this relation is independent of the parametrization.
		
		Further, given any parametrization for a family of lines as $C(t) = \{ \vc{x} = \lambda (\sigma,t), \sigma\in I\subset \mathbb{R}\}$, we can define the velocity of the line as
		\begin{equation}
			\label{contact_line_velocity}
			\vc{v}^c= \partial_t \lambda(\sigma,t).
		\end{equation}

		We also need the following two transport relations, on a domain $\Omega(t)$ changing with time, for any $w(\vc{x},t)$ smooth enough, defined on $\mathbb{R}^3$, and  
		$\varphi\in C^\infty_c(\mathbb{R}^3 \times \mathbb{R}_+)$,
		\begin{equation}
			\label{36a.}
			\int_0^\infty \int_{\Omega(t)}{w\; \partial_t \varphi \;\mathrm{d}t\;\mathrm{d}x} = - \int_0^\infty \int_{\Omega(t)}{\partial_t w \;\varphi \;\mathrm{d}t\;\mathrm{d}x} - \int_0^\infty \int_{\partial\Omega(t)}{w \; v^s_{n} \;\varphi \;\mathrm{d}t\; \mathrm{d}S}, 
		\end{equation}
		where $\vc{n}$ is the unit outward normal vector at the boundary $\partial\Omega(t)$, $v^s_n$ is the normal velocity of the boundary, and the last term on the RHS is a surface integral; Similarly, for any function $f$ defined over an interface $\Sigma(t)$, we have,
		\begin{align}
			\label{36b.}
			\int_0^\infty \int_{\Sigma(t)}{f\; \partial_t^\Sigma \varphi \;\mathrm{d}t\;\mathrm{d}S} & = - \int_0^\infty \int_{\Sigma(t)}{\partial_t^\Sigma f \varphi \;\mathrm{d}t\;\mathrm{d}S} - \int_0^\infty \int_{\partial\Sigma(t)}{f\; \vc{v}^c\cdot \nu \;\varphi \;\mathrm{d}t \;\mathrm{d}o}\nonumber\\
			& \qquad - \int_0^\infty \int_{\Sigma(t)}{f v^s_{n} \kappa_\Sigma\; \varphi\;\mathrm{d}t\;\mathrm{d}S},
		\end{align}
		where $\partial^\Sigma_t$ is the normal derivative over the surface, defined in (\ref{normal_derivative.}), $\kappa_\Sigma$ is the curvature of the surface, $\nu$ is the unit outward normal vector to the boundary $\partial \Sigma(t)$, and the second term on the RHS is a line integral.
		
		\subsection{Some explicit calculation}
		\label{explicit_terms}
		Here we write down explicitly several differential terms over the surface, in dimension two, for clarity.
		
		Let us consider an interface $y=h(x,t)$ with velocity $\vc{v}^s = (v^s_1,v^s_2)$. The unit outward normal vector at the surface is then given by
		\begin{equation}
			\label{normal}
			\vc{n} = \frac{\left( -\partial_x h,1\right) }{(1+ |\partial_xh|^2)^{1/2}}.
		\end{equation}
		
		The shape of the surface can be characterized in terms of its normal velocity (cf. (\ref{2})),
		\begin{equation}
			\label{normal_1}
			\partial_t h + v_1 \partial_x h -v_2=0.
		\end{equation}
		Indeed, at the free interface $f\equiv h(x,t)-y =0$, the normal is given by $\vc{n}_1 = \frac{\nabla_{\mathbf{x}} f}{|\nabla_{\mathbf{x}} f|}$ and hence, equation (\ref{2}) means precisely the convective derivative is zero, i.e.
		\begin{equation*}
			\partial_t f + \vc{v}^s \cdot \nabla_{\mathbf{x}} f=0.
		\end{equation*}
		
		In order to compute $\partial_t^\Sigma$, we choose a suitable normal parametrization of the surface,
		\begin{equation*}
			\gamma : \zeta\in\mathbb{R} \to \gamma(\zeta,t) \in \mathbb{R}^2 \quad \text{ such that } \quad \frac{\partial\gamma}{\partial t}\parallel\vc{n}.
		\end{equation*}
		We would like to relate $\gamma$ with the known parametrization $h$ (which is not necessarily a normal parametrization) and then express $\partial_t^\Sigma$ in terms of $h$. Let us introduce the following notations, for any generic function $w$ and a fixed point $(\overline{x}, \overline{y})$ on the interface $y = h(x,0)$,
		\begin{equation*}
			W(\overline{x},t) := w (\gamma(\overline{x},t),t) = w (\gamma_1(\overline{x},t),\gamma_2(\overline{x},t),t), \quad \quad \overline{y} = h(\overline{x},0)
		\end{equation*}
		and
		\begin{equation*}
			\tilde{w}(x,t) := w(x, h(x,t), t).
		\end{equation*}
		If we write the following Taylor expansion,
		\begin{equation}
			\label{54}
			\renewcommand{\arraystretch}{1.5}
			\gamma( \overline{x}, t )=
			\big(\gamma_1(\overline{x},t),\gamma_2(\overline{x},t)\big)
			=\big(\overline{x},\overline{y}\big) +  t\;v_n ( \overline{x}, 0) \vc{n}(\overline{x},t) + O (t^2),
		\end{equation}
		then
		\begin{equation*}
			\frac{\partial \gamma}{\partial t} = v_n ( \overline{x}, 0) \vc{n}(\overline{x},t) + O(t)
		\end{equation*}
		as desired. Further notice that
		\begin{equation*}
			\gamma_2(\overline{x},t) = h(\gamma_1(\overline{x},t)).
		\end{equation*}
		Indeed, from (\ref{54}), we have, combining (\ref{normal}) and (\ref{normal_1}),
		\begin{equation*}
			\begin{aligned}
				h(\gamma_1(\overline{x},t)) &= h\Big( \overline{x} - \frac{t}{(1+ |\partial_x h|^2)^{1/2}}\partial_x h(\overline{x},0)v_n(\overline{x},0)\Big) \\
				&= h(\overline{x},0) - \partial_xh(\overline{x},0) \frac{t}{(1+ |\partial_x h|^2)^{1/2}}\;\partial_x h(\overline{x},0)v_n(\overline{x},0)+t\; \partial_t h(\overline{x},0)\\\
				&= \overline{y} -\frac{t}{(1+ |\partial_x h|^2)^{1/2}} \; (\partial_x h(\overline{x},0))^2v_n(\overline{x},0)+ (1+ |\partial_x h|^2)^{1/2} \;v_n (\overline{x},0)\; t\\
				&= \overline{y} + \frac{t}{(1 + |\partial_x h|^2)^{1/2}}v_n(\overline{x},0) = \gamma_2(\overline{x},t).
			\end{aligned}
		\end{equation*}
		Therefore, we obtain
		\begin{equation*}
			W(\overline{x},t) = \tilde{w} \big( \gamma_1(\overline{x},t), t\big) =\tilde{w}\Big( \overline{x} - \frac{t}{(1+ |\partial_x h|^2)^{1/2}}\partial_x h(\overline{x},0)v_n(\overline{x},0),\; t\Big),
		\end{equation*}
		and then by definition,
		\begin{equation}
			\label{normal_derivative}
			\partial_t^\Sigma w = \frac{\mathrm{d} W}{\mathrm{d}t}\Big\rvert_{t=0}
			= -\frac{\partial_x h(\overline{x},0) v_n(\overline{x},0)}{(1+ |\partial_x h(\overline{x},0)|^2)^{1/2}} \frac{\partial \tilde{w}}{\partial x} + \frac{\partial \tilde{w}}{\partial t}.
		\end{equation}
		
		Next, for any tangential vector over the surface $A = a(s) \tau = a(x, h(x,t)) \vc{\tau}$ where $s$ is the arc length, defined by,
		\begin{equation*}
			s = \int_0^{x} (1+ {|\partial_x h(\zeta)|)^{1/2} \;\mathrm{d}\zeta }, \qquad \text{ hence, } \quad \mathrm{d}s = (1+ |\partial_x h|^2)^{1/2} \;\mathrm{d}x,
		\end{equation*}
		the surface divergence reduces to,
		\begin{equation*}
			\int_S {\mathrm{div}_\Sigma A \;\varphi \; \mathrm{d}s} = -\int_S {A \cdot \nabla_{\tau}\varphi \; \mathrm{d}s} = -\int_S {a \tau \cdot \frac{\mathrm{d}\varphi}{\mathrm{d}s}\tau\; \mathrm{d}s} = -\int_S {a \frac{\mathrm{d}\varphi}{\mathrm{d}s}\; \mathrm{d}s} 
		\end{equation*}
		which gives,
		\begin{equation}
			\label{surface_divergence}
			\mathrm{div}_\Sigma A = \frac{\mathrm{d}a }{\mathrm{d}s } = \frac{1}{(1+ |\partial_x h|^2)^{1/2}} \;\frac{\mathrm{d}}{\mathrm{d} x} \left( \tilde{a}(x)\right) 
		\end{equation}
		where we have denoted $\tilde{a}(x):= a(s)$.
		
		We can further obtain from (\ref{surface_divergence}), the following formula for tangential gradient, for a function $f$ on the surface $y=h(x,t)$,
		\begin{equation}
			\label{tan_grad}
			\nabla_{\tau}f \cdot\tau 
			= \frac{1}{(1+ |\partial_x h|^2)^{1/2}}\;\frac{\mathrm{d} \tilde{f}}{\mathrm{d} x}  ,
		\end{equation}
		where $\tilde{f}(x) := f(x, h(x,t))$. One can also obtain the above formula by considering any extension $\overline{f}$ in $\mathbb{R}^2$ and then computing $ \nabla \;\overline{f} \cdot\tau$, since the final result does not depend on the choice of extension.			
		
		\section{List of variables}
		
		The superscripts $\xi, s$ denote the quantities in the bulk and at the interfaces respectively, for
		$\xi= L,S,G,\ i=1,2,3$.
		\begin{align*}
			\Omega_\xi \qquad &\text{ region occupied by the bulk}\\
			\Sigma_i \qquad &\text{ interfaces between two bulk region}\\
			C \qquad &\text{ contact line}\\
			\chi^\xi \qquad &\text{ characteristic functions of the region occupied by the bulk}\\
			\delta_{\Sigma_i} \qquad &\text{ Dirac measures supported on the surfaces $\Sigma_i$}\\
			\varrho^\xi, \varrho^s_i, \qquad &\text{ mass densities}\\
			\vc{v}^\xi, \vc{v}^s_i, \qquad &\text{ velocities}\\
			\vc{v}^c \qquad &\text{ velocity of the contact line}\\
			\mathbb{S}^\xi, \mathbb{S}^s_i \qquad &\text{ stress tensors}\\
			e^\xi, e^s_i \qquad &\text{ total energy per unit mass}\\
			u^\xi, u^s_i \qquad &\text{ internal energy per unit mass}\\
			\vc{J}^\xi_q, \vc{J}^s_{q,i} \qquad &\text{ heat fluxes}\\
			\vc{n}_i, i=1,2 \qquad & \text{ unit normal vectors at $\Sigma_i$ inward with respect to the liquid domain}\\
			\vc{n}_3 \qquad & \text{ unit normal vector at $\Sigma_3$ inward with respect to the gas domain} \\
			\vc{w}_\tau, w_n \qquad &\text{ tangential and normal component of a vector $\vc{w}$ at each interface $\Sigma_i$} \\
			\vc{\nu}_i \qquad &\text{ vectors orthogonal to the contact line, tangent to the interfaces}\\
			s^\xi, s^s_i \qquad &\text{ specific entropies per unit mass}\\
			\vc{J}_s^\xi, \vc{J}^s_{s,i} \qquad &\text{ entropy fluxes per unit mass}\\
			\sigma^\xi, \sigma^s_i \qquad &\text{ entropy production per unit volume }\\
			T \qquad &\text{ absolute temperature}\\
			\mu^\xi, \mu^s_i \qquad &\text{ chemical potentials per unit mass}\\
			P^\xi, P^s_i \qquad &\text{ mechanical pressure}\\
			p^\xi \qquad &\text{ thermodynamic pressure}\\
			p^s_i \qquad &\text{ surface tension}\\
			\partial^\Sigma_t \qquad &\text{ normal derivative over the interface}\\
			\mathrm{div}_\Sigma, \nabla_{\tau} \qquad &\text{ surface divergence and tangential gradient}\\
			\kappa_\Sigma \qquad &\text{ curvature of the interface}\\
			\mu \qquad &\text{ dynamic viscosity of the fluid}\\
			\vc{V}_n \qquad &\text{ normal velocity of the liquid-solid interface}\\
			\tau_i \qquad &\text{ relaxation time of the interface $\Sigma_i$}\\
			\varrho^s_{i,e} \qquad &\text{ equilibrium surface density}\\
			\varrho^s_{i,0} \qquad &\text{ surface density corresponding to zero surface pressure}\\
			\theta \qquad &\text{ contact angle}\\
			\begin{rcases*}
				a^\xi, a^s_i, b^\xi, \\
				b^s_i, c^\xi, c^s_i, \\
				\alpha_i, \beta_i, \gamma_i 
			\end{rcases*}
			\quad &\text{ phenomenological constants at the interfaces}.
		\end{align*}
		
		\bigskip
		
		\textbf{Acknowledgements.} The authors acknowledge the support of the Hausdorff Center of Mathematics at the University of Bonn, funded by the Deutsche Forschungsgemeinschaft (DFG, German Research
		Foundation) under Germany's Excellence Strategy – EXC-2047/1 – 390685813.
		
		The authors certify that they do not have any conflict of interest.
		
		{\small

		}
	\end{document}